# FIXED-DOMAIN ASYMPTOTICS FOR A SUBCLASS OF MATÉRN-TYPE GAUSSIAN RANDOM FIELDS

By Wei-Liem Loh

*National University of Singapore*

Stein [*Statist. Sci.* **4** (1989) 432–433] proposed the Matérn-type Gaussian random fields as a very flexible class of models for computer experiments. This article considers a subclass of these models that are exactly once mean square differentiable. In particular, the likelihood function is determined in closed form, and under mild conditions the sieve maximum likelihood estimators for the parameters of the covariance function are shown to be weakly consistent with respect to fixed-domain asymptotics.

**1. Introduction.** In the modeling of computer experiments, it has become rather common practice to approximate the deterministic response as a realization of a stochastic process. In this regard, Sacks, Welch, Mitchell and Wynn [10] and Sacks, Schiller and Welch [9] proposed modeling using a Gaussian random field $X(x)$, $x \in [0,1]^d$, with a multiplicative covariance function,

$$\operatorname{Cov}(X(x), X(y)) = \sigma^2 \prod_{t=1}^{d} \exp(-\theta_t |x_t - y_t|^\gamma) \qquad \forall\, x = (x_1, \ldots, x_d)', y = (y_1, \ldots, y_d)' \in [0,1]^d, \tag{1}$$

where $\gamma \in (0,2], \theta_1, \ldots, \theta_d$ and $\sigma^2$ are strictly positive parameters. Ying [16, 17] investigated the fixed-domain asymptotic properties of the maximum likelihood estimators of the covariance function when $\gamma = 1$. In particular, he proved that the estimators are strongly consistent and asymptotically normal under mild conditions. van der Vaart [13] showed that when $\gamma = 1$ and $d = 2$, the maximum likelihood estimators are also asymptotically efficient. Recently, Loh and Lam [6] showed that, when $\gamma = 2$, sieve maximum likelihood estimators of $\theta_1, \ldots, \theta_d$ are strongly consistent using a regular sampling grid.

Stein [11] observed that the Gaussian model given by (1) may have some undesirable properties. In particular, for $\gamma \in (0, 2)$, the Gaussian random field with covariance function as in (1) will not be mean square differentiable. However, for the case $\gamma = 2$, it is infinitely mean square differentiable. Not allowing for processes that are differentiable but not infinitely differentiable may be unnecessarily restrictive. Stein then suggested using a Gaussian random field model, $X(x)$, $x \in [0,1]^d$, with the multiplicative Matérn-type covariance function

$$\begin{aligned}\operatorname{Cov}(X(x), X(y)) = \prod_{t=1}^{d} \frac{\pi^{1/2}\phi}{2^{\alpha-1}\Gamma(\alpha+1/2)\theta_t^{2\alpha}} (\theta_t|x_t - y_t|)^{\alpha} K_{\alpha}(\theta_t|x_t - y_t|) \\ \forall\, x = (x_1, \ldots, x_d)', y = (y_1, \ldots, y_d)' \in [0,1]^d,\end{aligned} \tag{2}$$

where $\alpha, \phi, \theta_1, \ldots, \theta_d$ are positive constants and $K_\alpha$ is the modified Bessel function of the second kind (see, e.g., [3], pages 222–223). The interesting parameter is $\alpha$, where $X$ will be $m$ times mean square differentiable if and only if $\alpha > m$.

Due to the dependence of the models with a multiplicative covariance function on the choice of coordinate axes, Stein ([12], pages 48–55) later advocated the use of Gaussian random field models with isotropic Matérn-type covariance functions. This class of covariance functions was first proposed by Matérn in 1960 as a reasonable class of models for isotropic random fields (see [7]). Stein ([12], Section 6.7) investigated the performance of maximum likelihood estimators for the parameters of a periodic version of the Matérn model with the hope that the large sample results for this periodic model will be similar to those for nonperiodic Matérn-type Gaussian random fields under fixed-domain asymptotics.

This article assumes that $X(x), x \in [0,1]^d$, is a mean-zero Gaussian random field with covariance function given as in (2), where $\alpha = 3/2$ and $\phi, \theta_1, \ldots, \theta_d$ are unknown positive constants. Since $K_{3/2}(s) = (\pi s^{-3}/2)^{1/2}(1 + s)e^{-s}$ (see [15], page 747), we observe that the covariance function of $X(x)$ is

$$\begin{aligned}\operatorname{Cov}(X(x), X(y)) = \frac{\pi^d \phi^d}{2^d \theta_1^3 \cdots \theta_d^3} \prod_{t=1}^{d} (1 + \theta_t|x_t - y_t|) e^{-\theta_t|x_t - y_t|} \\ \forall\, x = (x_1, \ldots, x_d)', y = (y_1, \ldots, y_d)' \in [0,1]^d.\end{aligned} \tag{3}$$

REMARK. For reasons of mathematical tractability, we have assumed a multiplicative Matérn-type covariance function and not an isotropic one.

The class of multiplicative Matérn-type covariance functions has nonetheless been used in, for example, [14]. A referee noted that the multiplicative Matérn-type covariance function and an isotropic one have a major difference under fixed-domain asymptotics. In particular, Zhang [18] recently proved that consistent estimators do not exist for all parameters in the latter case for $d = 1, 2$ or 3.

We are concerned with the estimation of $\phi, \theta_1, \ldots, \theta_d$ using observations that are taken from the above random field on a regular grid, that is,

$$\left\{X\left(\frac{i_1}{n}, \ldots, \frac{i_d}{n}\right) : 1 \le i_t \le n, 1 \le t \le d\right\}, \tag{4}$$

where $n$ is a strictly positive integer.

For simplicity, we order the elements of the set (4) lexicographically as an $n^d \times 1$ column vector $\tilde{X}_n$. Thus, the element $X(i_1/n, \ldots, i_d/n)$ precedes the element $X(j_1/n, \ldots, j_d/n)$ in $\tilde{X}_n$ if and only if there exists a $1 \le k \le d$ such that $i_t = j_t$ whenever $1 \le t < k$ and $i_k < j_k$. Then the covariance matrix $\Sigma_{\phi,\theta_1,\ldots,\theta_d;n}$ of $\tilde{X}_n$ is given by

$$\Sigma_{\phi,\theta_1,\ldots,\theta_d;n} = \frac{\pi^d \phi^d}{2^d \theta_1^3 \cdots \theta_d^3} \bigotimes_{t=1}^{d} R_{\theta_t,n},$$

where the symbol "$\bigotimes$" denotes the Kronecker product (see [2], page 599), and for each $1 \le t \le d$ $R_{\theta_t,n}$ denotes the $n \times n$ matrix whose $(i,j)$th element is $(1 + \theta_t |i-j|/n) \exp[-\theta_t |i-j|/n]$. Since $\tilde{X}_n \sim N_{n^d}(0, \Sigma_{\phi,\theta_1,\ldots,\theta_d;n})$, the likelihood function is

$$L_n(\phi, \theta_1, \ldots, \theta_d) = (2\pi)^{-n^d/2} |\Sigma_{\phi,\theta_1,\ldots,\theta_d;n}|^{-1/2} \exp(-\tilde{X}_n' \Sigma^{-1}_{\phi,\theta_1,\ldots,\theta_d;n} \tilde{X}_n / 2),$$

and the log-likelihood satisfies

$$\begin{aligned} &2 \log L_n(\phi, \theta_1, \ldots, \theta_d) \\ &\quad = -n^d \log(2\pi) - n^d \log\left(\frac{\pi^d \phi^d}{2^d \theta_1^3 \cdots \theta_d^3}\right) \\ &\qquad - \log\left|\bigotimes_{t=1}^{d} R_{\theta_t,n}\right| - \frac{2^d \theta_1^3 \cdots \theta_d^3}{\pi^d \phi^d} \tilde{X}_n' \left(\bigotimes_{t=1}^{d} R_{\theta_t,n}\right)^{-1} \tilde{X}_n. \end{aligned} \tag{5}$$

The rest of this article is organized as follows. Theorem 1 in Section 2 gives an exact closed-form formula for the determinant of the matrix $R_{\theta_t,n}$. Similarly, Theorems 2–4 in Section 3 establish an exact closed-form formula for the inverse of $R_{\theta_t,n}$. It should be noted that the closed-form formulas (especially the inverse) are not simple but fortunately are amenable to theoretical analysis with the help of a mathematical software system such as *Mathematica* [15] that has symbolic computation capability.

In Section 4 we simplify the exact formulas of the previous two sections via asymptotic approximations. These approximations are very sharp in that the error is of the order $O(c^n)$ for some constant $0 < c < 1$. In particular, Theorem 5 gives an asymptotic approximation of the determinant $|R_{\theta_t,n}|$ and Theorems 6 and 7 give asymptotic approximations to the elements of the inverse $R_{\theta_t,n}^{-1}$.

Using the results of Section 4 and *Mathematica*, Theorem 8 in Section 5 computes a large sample approximation of the Fisher information matrix for the parameters $\phi, \theta_1, \ldots, \theta_d$.

Section 6 gives a simple consistent maximum likelihood-type estimate $\hat{\phi}_{\tilde{\theta}_1,\ldots,\tilde{\theta}_d}$ (see Theorem 9) for $\phi$ even if the parameters $\theta_1, \ldots, \theta_d$ are misspecified. Section 7 defines a sieve maximum likehood estimator $(\hat{\phi}, \hat{\theta}_1, \ldots, \hat{\theta}_d)$ for $(\phi, \theta_1, \ldots, \theta_d)$. Theorem 10 establishes the weak consistency of $(\hat{\phi}, \hat{\theta}_1, \ldots, \hat{\theta}_d)$ under mild conditions when $d \geq 3$.

Appendix A contains technical results that are needed in the main body of this article. Many of the results found in Appendix A, though conceptually simple, involve extremely complicated computations and series expansions and the use of *Mathematica* is critical here. Appendix B contains the definitions of fourteen constants that are used in the exact expression of $R_{\theta_t,n}^{-1}$ in Section 3.

**2. Determinant.** In this section our main aim is to evaluate the determinant of $R_{\theta,n}$ for some strictly positive constant $\theta$. For simplicity, we write $w = \theta/n, u = e^{-w}$, and for $1 \leq m \leq n$ $R_{\theta,n;m}$ denotes the $m \times m$ matrix whose $(i,j)$th element is $(R_{\theta,n;m})_{i,j} = (1 + |i-j|w)u^{|i-j|}$, for all $1 \leq i, j \leq m$. Clearly, we have $R_{\theta,n;n} = R_{\theta,n}$.

PROPOSITION 1. *The determinant of $R_{\theta,n;m}$ is given by*

$$|R_{\theta,n;1}| = 1,$$

*and for $2 \leq m \leq n$ it satisfies the recurrence relation*

$$|R_{\theta,n;m}| = (-1)^{m-2}\{[1 + (m-2)w] - (1+w)[1 + (m-1)w]u^2\}u^{m-2}\tau_{-1}^{m-2} + \sum_{k=0}^{m-3} (-1)^k \tau_k |R_{\theta,n;m-k-1}| \tau_{-1}^k,$$

*where*

$$\tau_k = \begin{cases} 0, & \text{if } k = -2, \\ (w-1)u + (1+w)u^3, & \text{if } k = -1, \\ (1+kw)u^k - 2[1 + (k+1)w]u^{k+2} + [1 + (k+2)w]u^{k+4}, & \text{if } k \geq 0. \end{cases}$$

PROOF. Using linear algebra elementary row operations, we observe that, for each $1 \leq m \leq n$, $R_{\theta,n;m}$ can be reduced to the $m \times m$ almost upper triangular matrix $A_{\theta,n;m}$, where

$$A_{\theta,n;m} = \begin{pmatrix} 1 & (1+w)u & (1+2w)u^2 & (1+3w)u^3 & \cdots & (1+(m-1)w)u^{m-1} \\ (1+w)u & 1 & (1+w)u & (1+2w)u^2 & \cdots & (1+(m-2)w)u^{m-2} \\ 0 & \tau_{-1} & \tau_0 & \tau_1 & \cdots & \tau_{m-3} \\ 0 & 0 & \tau_{-1} & \tau_0 & \cdots & \tau_{m-4} \\ \vdots & \vdots & \vdots & \vdots & \ddots & \vdots \\ 0 & 0 & 0 & 0 & \cdots & \tau_1 \\ 0 & 0 & 0 & 0 & \cdots & \tau_0 \end{pmatrix}.$$

Now use elementary row operations again to reduce $A_{\theta,n;m}$ to a diagonal matrix and then use the fact that the determinant of a diagonal matrix is equal to the product of its diagonal elements. This proves Proposition 1. □

The following theorem gives a closed form expression for $|R_{\theta,n;m}|$.

THEOREM 1. *For* $2 \le m \le n$,

$$|R_{\theta,n;m}| = \frac{(\tilde{a}\alpha_1 + u^2\tilde{b})[-(b+\tau_{-1}u)\alpha_1]^{m-2} - (\tilde{a}\alpha_2 + u^2\tilde{b})[-(b+\tau_{-1}u)\alpha_2]^{m-2}}{\alpha_1 - \alpha_2},$$

*where*

$$\begin{aligned} a &= 1 - 2u^2 + u^4 - 2u^2w + 2u^4w, \\ b &= 2u^2 - u^4 + w - u^4w - 1, \\ \tilde{a} &= 1 - u^2 - 2u^2w - u^2w^2, \\ \tilde{b} &= u^2 + w + u^2w - 1, \\ \alpha_1 &= \frac{-(a - 2\tau_{-1}u) - \sqrt{(a-2\tau_{-1}u)^2 - 4\tau_{-1}u(b+\tau_{-1}u)}}{2(b+\tau_{-1}u)}, \\ \alpha_2 &= \frac{-(a - 2\tau_{-1}u) + \sqrt{(a-2\tau_{-1}u)^2 - 4\tau_{-1}u(b+\tau_{-1}u)}}{2(b+\tau_{-1}u)}. \end{aligned}$$

PROOF. It is convenient to define

$$\tau_k^* = \begin{cases} (1+kw)u^k - (1+w)[1+(k+1)w]u^{k+2}, & \text{if } k \ge 0, \\ 0, & \text{otherwise.} \end{cases} \tag{6}$$

Then $\tau_0^* = |R_{\theta,n;2}|$. We observe from Proposition 1 that, for $2 \le m \le n$,

$$|R_{\theta,n;m}| = \tau_{m-2}^*(-\tau_{-1})^{m-2} + \sum_{k=2}^{m-1} \tau_{m-k-1}(-\tau_{-1})^{m-k-1}|R_{\theta,n;k}|. \tag{7}$$

Applying (7) recursively to the right-hand side of itself, we obtain

$$
\begin{aligned}
|R_{\theta,n;m}| &= \tau^*_{m-2}(-\tau_{-1})^{m-2} + \tau_{m-3}\tau^*_0(-\tau_{-1})^{m-3} \\
&\quad + \sum_{k_1=3}^{m-1} \tau_{m-k_1-1}(-\tau_{-1})^{m-k_1-1}|R_{\theta,n;k_1}| \\
&= \tau^*_{m-2}(-\tau_{-1})^{m-2} + \tau_{m-3}\tau^*_0(-\tau_{-1})^{m-3} \\
&\quad + \sum_{k_1=3}^{m-1} \tau_{m-k_1-1}(-\tau_{-1})^{m-k_1-1} \\
&\qquad \times \Big\{ \tau^*_{k_1-2}(-\tau_{-1})^{k_1-2} + \tau_{k_1-3}\tau^*_0(-\tau_{-1})^{k_1-3} \\
&\qquad\quad + \sum_{k_2=3}^{k_1-1} \tau_{k_1-k_2-1}(-\tau_{-1})^{k_1-k_2-1}|R_{\theta,n;k_2}| \Big\} \\
&= (-\tau_{-1})^{m-2}\tau^*_{m-2} + (-\tau_{-1})^{m-3}\sum_{k_1=2}^{m-1}\tau_{m-k_1-1}\tau^*_{k_1-2} \\
&\quad + \sum_{k_1=3}^{m-1}\sum_{k_2=2}^{k_1-1}\tau_{m-k_1-1}\tau_{k_1-k_2-1}(-\tau_{-1})^{m-k_2-2}|R_{\theta,n;k_2}|.
\end{aligned}
$$

Now continuing this argument repeatedly, we obtain

$$
\begin{aligned}
|R_{\theta,n;m}| &= (-\tau_{-1})^{m-2}\tau^*_{m-2} + (-\tau_{-1})^{m-3}\sum_{k_1=2}^{m-1}\tau_{m-k_1-1}\tau^*_{k_1-2} \\
&\quad + (-\tau_{-1})^{m-4}\sum_{k_1=3}^{m-1}\sum_{k_2=2}^{k_1-1}\tau_{m-k_1-1}\tau_{k_1-k_2-1}\tau^*_{k_2-2} \\
&\quad + (-\tau_{-1})^{m-5}\sum_{k_1=4}^{m-1}\sum_{k_2=3}^{k_1-1}\sum_{k_3=2}^{k_2-1}\tau_{m-k_1-1}\tau_{k_1-k_2-1}\tau_{k_2-k_3-1}\tau^*_{k_3-2} \\
&\quad + \cdots + (-\tau_{-1})\sum_{k_1=m-2}^{m-1}\sum_{k_2=m-3}^{k_1-1}\cdots\sum_{k_{m-3}=2}^{k_{m-4}-1}\tau_{m-k_1-1}\tau_{k_1-k_2-1}\cdots \\
&\qquad \times \tau_{k_{m-4}-k_{m-3}-1}\tau^*_{k_{m-3}-2} \\
&\quad + \tau_0^{m-2}\tau^*_0.
\end{aligned}
\tag{8}
$$

Hence, it follows from (8) and Lemma A.1 (see Appendix A) that $|R_{\theta,n;m}|$ equals the coefficient of $z^{m-2}$ in the series expansion of

$$
\begin{aligned}
&(-\tau_{-1}u)^{m-2}G_{\tau^*}(z) + (-\tau_{-1}u)^{m-3}zG_\tau(z)G_{\tau^*}(z) + (-\tau_{-1}u)^{m-4}z^2G^2_\tau(z)G_{\tau^*}(z) \\
&\quad + \cdots + (-\tau_{-1}u)z^{m-3}G^{m-3}_\tau(z)G_{\tau^*}(z) + z^{m-2}G^{m-2}_\tau(z)G_{\tau^*}(z)
\end{aligned}
$$

$$
\begin{aligned}
&= \sum_{k=0}^{m-2} (-\tau_{-1}u)^{m-2-k} z^k G_\tau^k(z) G_{\tau^*}(z) \\
&= (-\tau_{-1}u)^{m-2} G_{\tau^*}(z) \frac{(-\tau_{-1}u)^{-m+1} z^{m-1} G_\tau^{m-1}(z) - 1}{(-\tau_{-1}u)^{-1} z G_\tau(z) - 1}.
\end{aligned}
$$

This implies that $|R_{\theta,n;m}|$ equals the coefficient of $z^{m-2}$ in the series expansion of

$$
\begin{aligned}
\frac{(-\tau_{-1}u)^{m-2} G_{\tau^*}(z)}{(\tau_{-1}u)^{-1} z G_\tau(z) + 1} &= \frac{(-\tau_{-1}u)^{m-2}(\tilde{a} + \tilde{b}z)}{(\tau_{-1}u)^{-1} z(a + bz) + (1 - z)^2} \\
&= \frac{(-1)^m (\tau_{-1}u)^{m-1}(\tilde{a} + \tilde{b}z)}{\tau_{-1}u + (a - 2\tau_{-1}u)z + (b + \tau_{-1}u)z^2}.
\end{aligned} \tag{9}
$$

Since $(a - 2\tau_{-1}u)^2 > 4\tau_{-1}u(b + \tau_{-1}u)$, the distinct roots of the quadratic equation

$$\tau_{-1}u + (a - 2\tau_{-1}u)z + (b + \tau_{-1}u)z^2 = 0$$

are $\alpha_1$ and $\alpha_2$. Hence, the right-hand side of (9) can be written as

$$
\begin{aligned}
&\frac{(-1)^m (\tau_{-1}u)^{m-1}(\tilde{a} + \tilde{b}z)}{(b + \tau_{-1}u)(\alpha_1 - z)(\alpha_2 - z)} \\
&\quad = \frac{(-1)^m (\tau_{-1}u)^{m-1}(\tilde{a} + \tilde{b}\alpha_1)}{(b + \tau_{-1}u)(\alpha_1 - z)(\alpha_2 - \alpha_1)} + \frac{(-1)^m (\tau_{-1}u)^{m-1}(\tilde{a} + \tilde{b}\alpha_2)}{(b + \tau_{-1}u)(\alpha_1 - \alpha_2)(\alpha_2 - z)} \\
&\quad = \frac{(-1)^m (\tau_{-1}u)^{m-1}(\tilde{a} + \tilde{b}\alpha_1)}{(b + \tau_{-1}u)\alpha_1(\alpha_2 - \alpha_1)} \sum_{k=0}^{\infty} \left(\frac{z}{\alpha_1}\right)^k \\
&\qquad + \frac{(-1)^m (\tau_{-1}u)^{m-1}(\tilde{a} + \tilde{b}\alpha_2)}{(b + \tau_{-1}u)(\alpha_1 - \alpha_2)\alpha_2} \sum_{k=0}^{\infty} \left(\frac{z}{\alpha_2}\right)^k.
\end{aligned}
$$

Thus, we conclude that, for $2 \le m \le n$,

$$
\begin{aligned}
|R_{\theta,n;m}| &= \frac{(-1)^m (\tau_{-1}u)^{m-1}(\tilde{a} + \tilde{b}\alpha_1)}{(b + \tau_{-1}u)\alpha_1^{m-1}(\alpha_2 - \alpha_1)} + \frac{(-1)^m (\tau_{-1}u)^{m-1}(\tilde{a} + \tilde{b}\alpha_2)}{(b + \tau_{-1}u)(\alpha_1 - \alpha_2)\alpha_2^{m-1}} \\
&= \frac{(-1)^m (b + \tau_{-1}u)^{m-2}[\tilde{a}(\alpha_2^{m-1} - \alpha_1^{m-1}) + \tilde{b}\alpha_1\alpha_2(\alpha_2^{m-2} - \alpha_1^{m-2})]}{\alpha_2 - \alpha_1},
\end{aligned}
$$

since $\alpha_1\alpha_2 = \tau_{-1}u/(b + \tau_{-1}u)$. This proves Theorem 1. □

**3. Inverse.** This section evaluates the inverse of $R_{\theta,n}$ for some positive constant $\theta$. For $1 \le i, j \le n$, if the $i$th row and $j$th column of $R_{\theta,n}$ are deleted, the resulting $(n-1) \times (n-1)$ matrix is denoted by $R_{\theta,n;-i,-j}$. Then $(-1)^{i+j}|R_{\theta,n;-i,-j}|$ is the cofactor of the $(i,j)$th element of the square matrix $R_{\theta,n}$. It is well known (see, e.g., [2], page 582) that

$$(R^{-1}_{\theta,n})_{i,j} = \frac{(-1)^{i+j}|R_{\theta,n;-i,-j}|}{|R_{\theta,n}|} \qquad \forall\, 1 \le i,j \le n. \tag{10}$$

Hence, it follows from Theorem 1 that it suffices to obtain a closed-form expression for $|R_{\theta,n;-i,-j}|$. We shall, without loss of generality, assume for the rest of this section that $i \le j$ and $i \ge n-j+1$ since $(R_{\theta,n})_{i,j} = (R_{\theta,n})_{n-j+1,n-i+1}$ and, hence, via symmetry, $(R^{-1}_{\theta,n})_{i,j} = (R^{-1}_{\theta,n})_{n-j+1,n-i+1}$. It is convenient to define

$$\hat{\tau}_k = \begin{cases} 2\tau_{-1}u, & \text{if } k=-2,\\ (1+w)u - 3(1+w)u^3 + 2(1+2w)u^5, & \text{if } k=-1,\\ (1+kw)u^k - 3(1+(k+2)w)u^{k+4} & \\ \qquad + 2(1+(k+3)w)u^{k+6}, & \text{if } k \ge 0, \end{cases}$$

$$\tilde{\tau}_k = \begin{cases} \tau_{-1}u/2, & \text{if } k=-2,\\ ((2w-1)u + (1+2w)u^5)/2, & \text{if } k=-1,\\ (2(1+kw)u^k - 3(1+(k+1)w)u^{k+2} & \\ \qquad + (1+(k+3)w)u^{k+6})/2, & \text{if } k \ge 0, \end{cases}$$

and for $1 \le m \le n-1$, the $m \times m$ matrix

$$R_{\theta,n;-i,-j;m} = ((R_{\theta,n;-i,-j})_{k,l})_{1 \le k,l \le m}.$$

Clearly, $R_{\theta,n;-i,-j;n-1} = R_{\theta,n;-i,-j}$.

Case 1. Suppose $i=1$. Using elementary row operations, we observe that, for each $1 \le m \le n-1$, $R_{\theta,n;-1,-j;m}$ can be reduced to the almost upper triangular $m \times m$ matrix $A_{\theta,n;-1,-j;m}$, where $(A_{\theta,n;-1,-j;m})_{k,l} = 0$ whenever $k \ge l+2$ and

$$(A_{\theta,n;-1,-j;m})_{k,l} = \begin{cases} (1+|l-k-1|w)u^{|l-k-1|}, & \text{if } 1 \le k \le 2 \wedge m, 1 \le l \le (j-1) \wedge m,\\ (1+(l-k)w)u^{l-k}, & \text{if } 1 \le k \le 2 \wedge m, j \le l \le m, \end{cases}$$

$$(A_{\theta,n;-1,-j;m})_{k,l} = \begin{cases} \tau_{l-k-1}, & \text{if } 3 \le k \le m, k-1 \le l \le (j-1) \wedge m,\\ \tau_{l-k}, & \text{if } 3 \le k \le m, j \le l \le m. \end{cases}$$

Case 2. Suppose $i=2$. Using elementary row operations, we observe that, for each $1 \le m \le n-1$, $R_{\theta,n;-2,-j;m}$ can be reduced to the almost upper triangular $m \times m$ matrix $A_{\theta,n;-2,-j;m}$, where $(A_{\theta,n;-2,-j;m})_{k,l} = 0$ whenever $k \ge l+2$ and

$$(A_{\theta,n;-2,-j;m})_{1,l} = \begin{cases} (1+|l-1|w)u^{|l-1|}, & \text{if } 1 \le l \le (j-1) \wedge m,\\ (1+lw)u^l, & \text{if } j \le l \le m, \end{cases}$$

$$(A_{\theta,n;-2,-j;m})_{2,l} = \begin{cases} (1+|l-3|w)u^{|l-3|}, & \text{if } 1 \le l \le (j-1)\wedge m, \\ (1+(l-2)w)u^{l-2}, & \text{if } j \le l \le m, \end{cases}$$

$$(A_{\theta,n;-2,-j;m})_{3,l} = \begin{cases} \tilde{\tau}_{l-4}, & \text{if } 2 \le l \le (j-1)\wedge m, \\ \tilde{\tau}_{l-3}, & \text{if } j \le l \le m, \end{cases}$$

$$(A_{\theta,n;-2,-j;m})_{k,l} = \begin{cases} \tau_{l-k-1}, & \text{if } 4 \le k \le m, k-1 \le l \le (j-1)\wedge m, \\ \tau_{l-k}, & \text{if } 4 \le k \le m, j \le l \le m. \end{cases}$$

Case 3. Suppose $3 \le i \le n$. Using elementary row operations, we observe that, for each $1 \le m \le n-1$, $R_{\theta,n;-i,-j;m}$ can be reduced to the almost upper triangular $m\times m$ matrix $A_{\theta,n;-i,-j;m}$, where $(A_{\theta,n;-i,-j;m})_{k,l} = 0$ whenever $k \ge l+2$ and

$$(A_{\theta,n;-i,-j;m})_{k,l} = \begin{cases} (1+|l-k|w)u^{|l-k|}, & \text{if } 1 \le k \le 2\wedge m, \\ & 1 \le l \le (j-1)\wedge m, \\ (1+(l-k+1)w)u^{l-k+1}, & \text{if } 1 \le k \le 2\wedge m, \\ & j \le l \le m, \end{cases}$$

$$(A_{\theta,n;-i,-j;m})_{k,l} = \begin{cases} \tau_{l-k}, & \text{if } 3 \le k \le (i-1)\wedge m, \\ & k-1 \le l \le (j-1)\wedge m, \\ \tau_{l-k+1}, & \text{if } 3 \le k \le (i-1)\wedge m, j \le l \le m, \end{cases}$$

$$(A_{\theta,n;-i,-j;m})_{i,l} = \begin{cases} \hat{\tau}_{l-i-1}, & \text{if } i-1 \le l \le (j-1)\wedge m, \\ \hat{\tau}_{l-i}, & \text{if } j \le l \le m, \end{cases}$$

$$(A_{\theta,n;-i,-j;m})_{i+1,l} = \begin{cases} \tilde{\tau}_{l-i-2}, & \text{if } i \le l \le (j-1)\wedge m, \\ \tilde{\tau}_{l-i-1}, & \text{if } j \le l \le m, \end{cases}$$

$$(A_{\theta,n;-i,-j;m})_{k,l} = \begin{cases} \tau_{l-k-1}, & \text{if } i+2 \le k \le m, k-1 \le l \le (j-1)\wedge m, \\ \tau_{l-k}, & \text{if } i+2 \le k \le m, j \le l \le m. \end{cases}$$

Proposition 2. *For $1 \le i \le j \le n, i \ge n-j+1$, the determinant of $R_{\theta,n;-i,-j;m}, 1 \le m \le n-1$, satisfies the recurrence relation*

$$\begin{aligned} &|R_{\theta,n;-i,-j;m}| \\ &\quad = \sum_{k=0}^{m-1} (-1)^k (A_{\theta,n;-i,-j;n-1})_{m-k,m} |R_{\theta,n;-i,-j;m-k-1}| \\ &\qquad \times \prod_{l=0}^{k-1} (A_{\theta,n;-i,-j;n-1})_{m-l,m-l-1}, \end{aligned}$$

*where* $|R_{\theta,n;-i,-j;0}| = 1$.

Proof. Since $|R_{\theta,n;-i,-j;k}| = |A_{\theta,n;-i,-j;k}|$ for all $1 \le k \le n-1$, Proposition 2 follows by using elementary row operations to reduce $A_{\theta,n;-i,-j;m}$ to a diagonal matrix. □

Theorems 2, 3 and 4 below compute exact expressions for the determinants $|R_{\theta,n;-i,-j}|$, where $i \le j$ and $i \ge n-j+1$.

THEOREM 2. *For* $2 \le m \le n-1$,

$$|R_{\theta,n;-1,-n;m}| = w^2 u^2 \tau_{-1}^{m-2}.$$

PROOF. We observe from Proposition 2 and Case 1 that, for $2 \le m \le n-1$,

$$\begin{aligned}|R_{\theta,n;-1,-n;m}| &= \sum_{k=0}^{m-1} (-1)^{m-k-1} (A_{\theta,n;-1,-n;n-1})_{k+1,m} |R_{\theta,n;-1,-n;k}| \\ &\qquad \times \prod_{l=0}^{m-k-2} (A_{\theta,n;-1,-n;n-1})_{m-l,m-l-1} \\ &= (A_{\theta,n;-1,-n;n-1})_{m,m} |R_{\theta,n;-1,-n;m-1}| \\ &= |R_{\theta,n;-1,-n;2}| \tau_{-1}^{m-2} \\ &= w^2 u^2 \tau_{-1}^{m-2}.\end{aligned}$$

This proves Theorem 2. $\square$

For $1 \le m \le n-1$, let $\tilde{A}_{\theta,n;m}$ denote the $m \times m$ matrix such that

$$\tilde{A}_{\theta,n;m} = \begin{pmatrix} \tau_0 & \tau_1 & \tau_2 & \cdots & \tau_{m-2} & \tau_{m-1} \\ \tau_{-1} & \tau_0 & \tau_1 & \cdots & \tau_{m-3} & \tau_{m-2} \\ 0 & \tau_{-1} & \tau_0 & \cdots & \tau_{m-2} & \tau_{m-3} \\ \vdots & \vdots & \vdots & \ddots & \vdots & \vdots \\ 0 & 0 & 0 & \cdots & \tau_0 & \tau_1 \\ 0 & 0 & 0 & \cdots & \tau_{-1} & \tau_0 \end{pmatrix}. \tag{11}$$

THEOREM 3. *For* $n-1 \le j \le n$,

$$|R_{\theta,n;-2,-j;n-1}| = \tau_{-1}^{j-4}\{\tilde{\tau}_{-2}[(1+2w)^2 u^4 - 1] + \tilde{\tau}_{-1}\tau_1^*\}|\tilde{A}_{\theta,n;n-j}|,$$

*where* $\tau_1^*$ *is as in* (6) *and* $|\tilde{A}_{\theta,n;0}| = 1$.

PROOF. We observe from Proposition 2 and Case 2 that for $n-1 \le j \le n$,

$$\begin{aligned}|R_{\theta,n;-2,-j;1}| &= 1, \\ |R_{\theta,n;-2,-j;2}| &= (1+w)u - (1+w)(1+2w)u^3, \\ |R_{\theta,n;-2,-j;3}| &= \sum_{k=0}^{2} (-1)^{2-k} (A_{\theta,n;-2,-j;n-1})_{k+1,3} |R_{\theta,n;-2,-j;k}|\end{aligned}$$

$$\times \prod_{l=0}^{1-k} (A_{\theta,n;-2,-j;n-1})_{3-l,2-l}$$

$$= (1+2w)^2 u^4 \tilde{\tau}_{-2} - \tilde{\tau}_{-2} + \tilde{\tau}_{-1}[(1+w)u - (1+w)(1+2w)u^3],$$

$$|R_{\theta,n;-2,-j;n-1}| = |R_{\theta,n;-2,-j;3}| \tau_{-1}^{j-4} |\tilde{A}_{\theta,n;n-j}|.$$

This proves Theorem 3. □

THEOREM 4. *With the notation of Appendix B, for $3 \le i \le n-1$,*

$$|R_{\theta,n;-i,-i,m}| = |R_{\theta,n;m}| \qquad \forall 1 \le m \le i-1,$$

*and*

$$|R_{\theta,n;-i,-i;i}| = \frac{C_{7,1}(\tilde{a}\alpha_1 + u^2\tilde{b})}{\alpha_1 - \alpha_2}[-(b+\tau_{-1}u)\alpha_1]^{i-3} - \frac{C_{7,2}(\tilde{a}\alpha_2 + u^2\tilde{b})}{\alpha_1 - \alpha_2}[-(b+\tau_{-1}u)\alpha_2]^{i-3}. \tag{12}$$

*For $4 \le i+1 \le m \le n-1$ we have*

$$\begin{aligned}
|R_{\theta,n;-i,-i;m}| = {} & \frac{C_{1,1}\tau_{-1}u^4(\tilde{a}\alpha_1 + u^2\tilde{b})}{(\alpha_1-\alpha_2)^2} \\
& \times \{(\alpha_1 - u^2)[-(b+\tau_{-1}u)\alpha_1]^{m-4} \\
& \quad - (\alpha_2 - u^2)[-(b+\tau_{-1}u)\alpha_1]^{i-3}[-(b+\tau_{-1}u)\alpha_2]^{m-i-1}\} \\
& - \frac{C_{1,2}\tau_{-1}u^4(\tilde{a}\alpha_2 + u^2\tilde{b})}{(\alpha_1-\alpha_2)^2} \\
& \times \{(\alpha_1 - u^2)[-(b+\tau_{-1}u)\alpha_1]^{m-i-1}[-(b+\tau_{-1}u)\alpha_2]^{i-3} \\
& \quad - (\alpha_2 - u^2)[-(b+\tau_{-1}u)\alpha_2]^{m-4}\} \\
& + \frac{C_{2,1}\tau_{-1}u^6(\tilde{a}\alpha_1 + u^2\tilde{b})}{(\alpha_1-\alpha_2)^2} \\
& \times \{[-(b+\tau_{-1}u)\alpha_1]^{m-4} \\
& \quad - [-(b+\tau_{-1}u)\alpha_1]^{i-3}[-(b+\tau_{-1}u)\alpha_2]^{m-i-1}\} \\
& - \frac{C_{2,2}\tau_{-1}u^6(\tilde{a}\alpha_2 + u^2\tilde{b})}{(\alpha_1-\alpha_2)^2} \\
& \times \{[-(b+\tau_{-1}u)\alpha_1]^{m-i-1}[-(b+\tau_{-1}u)\alpha_2]^{i-3} \\
& \quad - [-(b+\tau_{-1}u)\alpha_2]^{m-4}\}.
\end{aligned} \tag{13}$$

*Also, for $3 \le i \le n-1$,*

$$\begin{aligned}
|R_{\theta,n;-i,-(i+1);m}| &= |R_{\theta,n;m}| \qquad \forall 1 \le m \le i-1, \\
|R_{\theta,n;-i,-(i+1);i}| &= \frac{C_{3,1}(\tilde{a}\alpha_1 + u^2\tilde{b})}{\alpha_1 - \alpha_2}[-(b+\tau_{-1}u)\alpha_1]^{i-3} \\
& \quad - \frac{C_{3,2}(\tilde{a}\alpha_2 + u^2\tilde{b})}{\alpha_1 - \alpha_2}[-(b+\tau_{-1}u)\alpha_2]^{i-3}
\end{aligned} \tag{14}$$

*and for* $4 \le i+1 \le m \le n-1$, *we have*

$$
\begin{aligned}
&|R_{\theta,n;-i,-(i+1);m}| \\
&\quad = \frac{C_{4,1}\tau_{-1}u^4(\tilde{a}\alpha_1 + u^2\tilde{b})}{(\alpha_1-\alpha_2)^2} \\
&\qquad \times \{(\alpha_1 - u^2)[-(b+\tau_{-1}u)\alpha_1]^{m-4} \\
&\qquad\quad - (\alpha_2 - u^2)[-(b+\tau_{-1}u)\alpha_1]^{i-3}[-(b+\tau_{-1}u)\alpha_2]^{m-i-1}\} \\
&\qquad - \frac{C_{4,2}\tau_{-1}u^4(\tilde{a}\alpha_2 + u^2\tilde{b})}{(\alpha_1-\alpha_2)^2} \\
&\qquad\quad \times \{(\alpha_1 - u^2)[-(b+\tau_{-1}u)\alpha_1]^{m-i-1}[-(b+\tau_{-1}u)\alpha_2]^{i-3} \\
&\qquad\qquad - (\alpha_2 - u^2)[-(b+\tau_{-1}u)\alpha_2]^{m-4}\} \\
&\qquad + \frac{C_{5,1}\tau_{-1}u^6(\tilde{a}\alpha_1 + u^2\tilde{b})}{(\alpha_1-\alpha_2)^2} \\
&\qquad\quad \times \{[-(b+\tau_{-1}u)\alpha_1]^{m-4} \\
&\qquad\qquad - [-(b+\tau_{-1}u)\alpha_1]^{i-3}[-(b+\tau_{-1}u)\alpha_2]^{m-i-1}\} \\
&\qquad - \frac{C_{5,2}\tau_{-1}u^6(\tilde{a}\alpha_2 + u^2\tilde{b})}{(\alpha_1-\alpha_2)^2} \\
&\qquad\quad \times \{[-(b+\tau_{-1}u)\alpha_1]^{m-i-1}[-(b+\tau_{-1}u)\alpha_2]^{i-3} \\
&\qquad\qquad - [-(b+\tau_{-1}u)\alpha_2]^{m-4}\}.
\end{aligned} \tag{15}
$$

*Finally, for* $3 \le i \le n-2$,

$$
\begin{aligned}
|R_{\theta,n;-i,-n;n-1}| &= \frac{C_{6,1}\tau_{-1}^{n-i-2}(\tilde{a}\alpha_1 + u^2\tilde{b})}{\alpha_1-\alpha_2}[-(b+\tau_{-1}u)\alpha_1]^{i-3} \\
&\quad - \frac{C_{6,2}\tau_{-1}^{n-i-2}(\tilde{a}\alpha_2 + u^2\tilde{b})}{\alpha_1-\alpha_2}[-(b+\tau_{-1}u)\alpha_2]^{i-3}
\end{aligned}
$$

*and for* $5 \le i+2 \le j \le n-1$,

$$
\begin{aligned}
|R_{\theta,n;-i,-j;n-1}| &= \frac{C_{6,1}\tau_{-1}^{j-i-2}(\tilde{a}\alpha_1 + u^2\tilde{b})}{(\alpha_1-\alpha_2)^2} \\
&\quad \times \{(a\alpha_1 + u^2b)[-(b+\tau_{-1}u)\alpha_1]^{n-j+i-4} \\
&\qquad - (a\alpha_2 + u^2b)[-(b+\tau_{-1}u)\alpha_1]^{i-3}[-(b+\tau_{-1}u)\alpha_2]^{n-j-1}\} \\
&\quad - \frac{C_{6,2}\tau_{-1}^{j-i-2}(\tilde{a}\alpha_2 + u^2\tilde{b})}{(\alpha_1-\alpha_2)^2} \\
&\qquad \times \{(a\alpha_1 + u^2b)[-(b+\tau_{-1}u)\alpha_1]^{n-j-1} \\
&\qquad\quad \times [-(b+\tau_{-1}u)\alpha_2]^{i-3} \\
&\qquad\quad - (a\alpha_2 + u^2b)[-(b+\tau_{-1}u)\alpha_2]^{n-j+i-4}\}.
\end{aligned} \tag{16}
$$

PROOF. The first equality is immediate from the definition of $R_{\theta,n;-i,-i;m}$. To prove (12), we observe from Proposition 2, Case 3 and (33) that, for

$3 \le i \le n-1$,

$$\begin{aligned}
|R_{\theta,n;-i,-i;i}| &= \sum_{k=0}^{i-1}(-1)^{i-k-1}(A_{\theta,n;-i,-i;n-1})_{k+1,i}|R_{\theta,n;-i,-i;k}| \\
&\qquad \times \prod_{l=0}^{i-k-2}(A_{\theta,n;-i,-i;n-1})_{i-l,i-l-1} \\
&= (-1)^{i-2}\tau_{i-1}^{*}\tau_{-1}^{i-3}\hat{\tau}_{-2} \\
&\quad + \sum_{k=2}^{i-2}(-1)^{i-k-1}\tau_{i-k}\tau_{-1}^{i-k-2}\hat{\tau}_{-2}|R_{\theta,n;k}| + \hat{\tau}_0|R_{\theta,n;i-1}| \\
&= (-1)^{i-2}\tau_{i-1}^{*}\tau_{-1}^{i-3}\hat{\tau}_{-2} + \hat{\tau}_0|R_{\theta,n;i-1}| \\
&\quad + \frac{(-1)^{i-1}\tau_{-1}^{i-3}\hat{\tau}_{-2}u^{i-1}(1-u^2)^2(\tilde{a}\alpha_1+u^2\tilde{b})}{(\alpha_1-\alpha_2)[\alpha_1(b+\tau_{-1}u)-\tau_{-1}u]} \\
&\quad \times \left\{\left[\frac{\alpha_1(b+\tau_{-1}u)}{\tau_{-1}u}\right]^{i-3}-1\right\} \\
&\quad - \frac{(-1)^{i-1}\tau_{-1}^{i-3}\hat{\tau}_{-2}u^{i-1}(1-u^2)^2(\tilde{a}\alpha_2+u^2\tilde{b})}{(\alpha_1-\alpha_2)[\alpha_2(b+\tau_{-1}u)-\tau_{-1}u]} \\
&\quad \times \left\{\left[\frac{\alpha_2(b+\tau_{-1}u)}{\tau_{-1}u}\right]^{i-3}-1\right\} \\
&\quad + \frac{(-1)^{i-1}\tau_{-1}^{i-2}\hat{\tau}_{-2}[(\tilde{a}\alpha_1+u^2\tilde{b})S_{i,i,1}-(\tilde{a}\alpha_2+u^2\tilde{b})S_{i,i,2}]}{\alpha_1-\alpha_2},
\end{aligned}$$

where $S_{i,m,l}$ is as in (31) and (12) follows after some algebra. For $4 \le i+1 \le m \le n-1$,

$$\begin{aligned}
|R_{\theta,n;-i,-i;m}| &= \sum_{k_1=0}^{m-1}(-1)^{m-k_1-1}(A_{\theta,n;-i,-i;n-1})_{k_1+1,m}|R_{\theta,n;-i,-i;k_1}| \\
&\qquad \times \prod_{l=0}^{m-k_1-2}(A_{\theta,n;-i,-i;n-1})_{m-l,m-l-1} \\
&= (-\tau_{-1})^m\xi_m \\
&\quad + (-\tau_{-1})^{m-1}\sum_{k_1=i+1}^{m-1}\tau_{m-k_1-1}(-\tau_{-1})^{-k_1}|R_{\theta,n;-i,-i;k_1}|,
\end{aligned}$$

where

$$
\begin{aligned}
\xi_m &= \tau^*_{m-1}\tau_{-1}^{-4}\hat{\tau}_{-2}\tilde{\tau}_{-1} + \sum_{k=2}^{i-2}(-1)^{-k-1}\tau_{m-k}\tau_{-1}^{-k-3}\hat{\tau}_{-2}\tilde{\tau}_{-1}|R_{\theta,n;k}| \\
&\quad + (-1)^{-i}\hat{\tau}_{m-i}\tau_{-1}^{-i-1}\tilde{\tau}_{-1}|R_{\theta,n;i-1}| \\
&\quad + (-1)^{-i-1}\tilde{\tau}_{m-i-1}\tau_{-1}^{-i-1}|R_{\theta,n;-i,-i;i}| \\
&= \tau^*_{m-1}\tau_{-1}^{-4}\hat{\tau}_{-2}\tilde{\tau}_{-1} + (-1)^{-i}\hat{\tau}_{m-i}\tau_{-1}^{-i-1}\tilde{\tau}_{-1}|R_{\theta,n;i-1}| \\
&\quad + (-1)^{-i-1}\tilde{\tau}_{m-i-1}\tau_{-1}^{-i-1}|R_{\theta,n;-i,-i;i}| \\
&\quad - \frac{\hat{\tau}_{-2}\tilde{\tau}_{-1}u^{m-1}(1-u^2)^2(\tilde{a}\alpha_1+u^2\tilde{b})}{\tau_{-1}^4(\alpha_1-\alpha_2)[\alpha_1(b+\tau_{-1}u)-\tau_{-1}u]}\left\{\left[\frac{\alpha_1(b+\tau_{-1}u)}{\tau_{-1}u}\right]^{i-3}-1\right\} \\
&\quad + \frac{\hat{\tau}_{-2}\tilde{\tau}_{-1}u^{m-1}(1-u^2)^2(\tilde{a}\alpha_2+u^2\tilde{b})}{\tau_{-1}^4(\alpha_1-\alpha_2)[\alpha_2(b+\tau_{-1}u)-\tau_{-1}u]}\left\{\left[\frac{\alpha_2(b+\tau_{-1}u)}{\tau_{-1}u}\right]^{i-3}-1\right\} \\
&\quad - \frac{\hat{\tau}_{-2}\tilde{\tau}_{-1}[(\tilde{a}\alpha_1+u^2\tilde{b})S_{i,m,1}-(\tilde{a}\alpha_2+u^2\tilde{b})S_{i,m,2}]}{\tau_{-1}^3(\alpha_1-\alpha_2)}.
\end{aligned}
\tag{17}
$$

The last equality uses (33). An immediate consequence is $|R_{\theta,n;-i,-i;i+1}| = (-\tau_{-1})^{i+1}\xi_{i+1}$. Then repeating the above argument, we obtain, for $4 \le i + 1 < m \le n-1$,

$$
\begin{aligned}
|R_{\theta,n;-i,-i;m}| &= (-\tau_{-1})^m\xi_m \\
&\quad + (-\tau_{-1})^{m-1}\sum_{k_1=i+1}^{m-1}\tau_{m-k_1-1}(-\tau_{-1})^{-k_1} \\
&\quad\quad \times \Bigg[(-\tau_{-1})^{k_1}\xi_{k_1} \\
&\quad\quad\quad + (-\tau_{-1})^{k_1-1}\sum_{k_2=i+1}^{k_1-1}\tau_{k_1-k_2-1}(-\tau_{-1})^{-k_2}|R_{\theta,n;-i,-i;k_2}|\Bigg] \\
&= (-\tau_{-1})^m\xi_m + (-\tau_{-1})^{m-1}\sum_{k_1=i+1}^{m-1}\tau_{m-k_1-1}\xi_{k_1} \\
&\quad + (-\tau_{-1})^{m-2}\sum_{k_1=i+2}^{m-1}\sum_{k_2=i+1}^{k_1-1}\tau_{m-k_1-1}\tau_{k_1-k_2-1} \\
&\quad\quad \times (-\tau_{-1})^{-k_2}|R_{\theta,n;-i,-i;k_2}| \\
&= (-\tau_{-1})^m\xi_m + (-\tau_{-1})^{m-1}\sum_{k_1=i+1}^{m-1}\tau_{m-k_1-1}\xi_{k_1} \\
&\quad + (-\tau_{-1})^{m-2}\sum_{k_1=i+2}^{m-1}\sum_{k_2=i+1}^{k_1-1}\tau_{m-k_1-1}\tau_{k_1-k_2-1}\xi_{k_2}
\end{aligned}
$$

$$
\begin{aligned}
&+ (-\tau_{-1})^{m-3} \sum_{k_1=i+3}^{m-1} \sum_{k_2=i+2}^{k_1-1} \sum_{k_3=i+1}^{k_2-1} \tau_{m-k_1-1}\tau_{k_1-k_2-1} \\
&\qquad\qquad \times \tau_{k_2-k_3-1}\xi_{k_3} + \cdots \\
&+ (-\tau_{-1})^{i+2} \\
&\quad \times \sum_{k_1=m-2}^{m-1} \sum_{k_2=m-3}^{k_1-1} \cdots \sum_{k_{m-i-2}=i+1}^{k_{m-i-3}-1} \tau_{m-k_1-1}\tau_{k_1-k_2-1}\cdots \\
&\qquad\qquad \times \tau_{k_{m-i-3}-k_{m-i-2}-1}\xi_{k_{m-i-2}} \\
&+ (-\tau_{-1})^{i+1}\tau_0^{m-i-1}\xi_{i+1}.
\end{aligned}
$$

Hence, it follows from Lemma A.4 (see Appendix A) that $|R_{\theta,n;-i,-i;m}|, 4 \le i+1 \le m \le n-1$, equals the coefficient of $z^m$ in the series expansion of

$$
\begin{aligned}
&(-\tau_{-1}u)^m G_\xi(z) + (-\tau_{-1}u)^{m-1} z G_\tau(z) G_\xi(z) \\
&\quad + (-\tau_{-1}u)^{m-2} z^2 G_\tau^2(z) G_\xi(z) + \cdots + (-\tau_{-1}u)^{i+1} z^{m-i-1} G_\tau^{m-i-1}(z) G_\xi(z) \\
&\quad = \sum_{k=0}^{m-i-1} (-\tau_{-1}u)^{m-k} z^k G_\tau^k(z) G_\xi(z) \\
&\quad = (-\tau_{-1}u)^m G_\xi(z) \sum_{k=0}^{m-i-1} (-\tau_{-1}u)^{-k} z^k G_\tau^k(z) \\
&\quad = (-\tau_{-1}u)^m G_\xi(z) \frac{(-\tau_{-1}u)^{-m+i} z^{m-i} G_\tau^{m-i}(z) - 1}{(-\tau_{-1}u)^{-1} z G_\tau(z) - 1}.
\end{aligned}
$$

It follows from Lemmas A.1 and A.4 that $|R_{\theta,n;-i,-i;m}|, 4 \le i+1 \le m \le n-1$, equals the coefficient of $z^m$ in the series expansion of

$$
\begin{aligned}
\frac{(-\tau_{-1}u)^m G_\xi(z)}{(\tau_{-1}u)^{-1} z G_\tau(z) + 1} &= \frac{z^{i+1}(-\tau_{-1}u)^m (a_\xi + b_\xi z)}{(\tau_{-1}u)^{-1} z(a+bz) + (1-z)^2} \\
&= \frac{(-1)^m z^{i+1} (\tau_{-1}u)^{m+1} (a_\xi + b_\xi z)}{\tau_{-1}u + (a - 2\tau_{-1}u) z + (b + \tau_{-1}u) z^2}.
\end{aligned}
\tag{18}
$$

Letting $\alpha_1$ and $\alpha_2$ be as in Theorem 1, the right-hand side of (18) can be written as

$$
\begin{aligned}
&\frac{(-1)^m z^{i+1} (\tau_{-1}u)^{m+1} (a_\xi + b_\xi z)}{(b + \tau_{-1}u)(\alpha_1 - z)(\alpha_2 - z)} \\
&\quad = \frac{(-1)^m (\tau_{-1}u)^{m+1} (a_\xi + b_\xi \alpha_1)}{(b + \tau_{-1}u)\alpha_1(\alpha_2 - \alpha_1)} \sum_{k=0}^{\infty} \frac{z^{k+i+1}}{\alpha_1^k}
\end{aligned}
$$

$$+\frac{(-1)^m(\tau_{-1}u)^{m+1}(a_\xi+b_\xi\alpha_2)}{(b+\tau_{-1}u)(\alpha_1-\alpha_2)\alpha_2}\sum_{k=0}^{\infty}\frac{z^{k+i+1}}{\alpha_2^k}.$$

Thus, we conclude that, for $4\le i+1\le m\le n-1$,

$$\begin{aligned}|R_{\theta,n;-i,-i,m}| &= \frac{(-1)^m(\tau_{-1}u)^{m+1}(a_\xi+b_\xi\alpha_1)}{(b+\tau_{-1}u)\alpha_1^{m-i}(\alpha_2-\alpha_1)}\\ &\quad+\frac{(-1)^m(\tau_{-1}u)^{m+1}(a_\xi+b_\xi\alpha_2)}{(b+\tau_{-1}u)(\alpha_1-\alpha_2)\alpha_2^{m-i}}\\ &= \frac{(-1)^m(\tau_{-1}u)^{i+1}a_\xi}{\alpha_1-\alpha_2}\{(\alpha_1-u^2)[(b+\tau_{-1}u)\alpha_1]^{m-i-1}\\ &\qquad -(\alpha_2-u^2)[(b+\tau_{-1}u)\alpha_2]^{m-i-1}\}\\ &\quad+\frac{(-1)^m(\tau_{-1}u)^{i+1}u^2(a_\xi+b_\xi)}{\alpha_1-\alpha_2}\\ &\qquad\times\{[(b+\tau_{-1}u)\alpha_1]^{m-i-1}-[(b+\tau_{-1}u)\alpha_2]^{m-i-1}\}.\end{aligned}\tag{19}$$

Equation (13) follows from (19) and Lemma A.4 after some algebra. Next, from the definition of $R_{\theta,n;-i,-(i+1);m}$ we have

$$|R_{\theta,n;-i,-(i+1);m}|=|R_{\theta,n;m}|\qquad\forall\,1\le m\le i-1,$$

and

$$\begin{aligned}&|R_{\theta,n;-i,-(i+1);i}|\\ &\quad=\sum_{k=0}^{i-1}(-1)^{i-k-1}(A_{\theta,n;-i,-(i+1);n-1})_{k+1,i}\\ &\qquad\qquad\times|R_{\theta,n;-i,-(i+1);k}|\prod_{l=0}^{i-k-2}(A_{\theta,n;-i,-(i+1);n-1})_{i-l,i-l-1}\\ &\quad=(-1)^{i-1}(1+w)[1+(i-1)w]u^i\tau_{-1}^{i-3}\hat{\tau}_{-2}\\ &\qquad+(-1)^{i-2}[1+(i-2)w]u^{i-2}\tau_{-1}^{i-3}\hat{\tau}_{-2}\\ &\qquad+\sum_{k=2}^{i-2}(-1)^{i-k-1}\tau_{i-k-1}\tau_{-1}^{i-k-2}\hat{\tau}_{-2}|R_{\theta,n;-i,-(i+1);k}|\\ &\qquad+\hat{\tau}_{-1}|R_{\theta,n;-i,-(i+1);i-1}|\\ &\quad=(-1)^{i-2}\tau_{i-2}^{*}\tau_{-1}^{i-3}\hat{\tau}_{-2}+\hat{\tau}_{-1}|R_{\theta,n;i-1}|\\ &\qquad+\sum_{k=2}^{i-2}(-1)^{i-k-1}\tau_{i-k-1}\tau_{-1}^{i-k-2}\hat{\tau}_{-2}|R_{\theta,n;k}|\\ &\quad=(-1)^{i-2}\tau_{i-2}^{*}\tau_{-1}^{i-3}\hat{\tau}_{-2}+\hat{\tau}_{-1}|R_{\theta,n;i-1}|\end{aligned}$$

$$
\begin{aligned}
&+\frac{(-1)^{i-1}\tau_{-1}^{i-2}\hat{\tau}_{-2}u^{i-2}(1-u^2)^2(\tilde{a}\alpha_1+u^2\tilde{b})}{\tau_{-1}(\alpha_1-\alpha_2)[(b+\tau_{-1}u)\alpha_1-\tau_{-1}u]}\left\{\left[\frac{(b+\tau_{-1}u)\alpha_1}{\tau_{-1}u}\right]^{i-3}-1\right\}\\
&-\frac{(-1)^{i-1}\tau_{-1}^{i-2}\hat{\tau}_{-2}u^{i-2}(1-u^2)^2(\tilde{a}\alpha_2+u^2\tilde{b})}{\tau_{-1}(\alpha_1-\alpha_2)[(b+\tau_{-1}u)\alpha_2-\tau_{-1}u]}\left\{\left[\frac{(b+\tau_{-1}u)\alpha_2}{\tau_{-1}u}\right]^{i-3}-1\right\}\\
&+\frac{(-1)^{i-1}\tau_{-1}^{i-2}\hat{\tau}_{-2}[(\tilde{a}\alpha_1+u^2\tilde{b})S_{i,i-1,1}-(\tilde{a}\alpha_2+u^2\tilde{b})S_{i,i-1,2}]}{\alpha_1-\alpha_2}.
\end{aligned}
$$

The last equality uses (33) and (14) results after some algebra. For $4 \le i+1 \le m \le n-1$, we have

$$
\begin{aligned}
|R_{\theta,n;-i,-(i+1);m}| &= \sum_{k=0}^{m-1}(-1)^{m-k-1}(A_{\theta,n;-i,-(i+1);n-1})_{k+1,m}\\
&\qquad \times |R_{\theta,n;-i,-(i+1);k}|\prod_{l=0}^{m-k-2}(A_{\theta,n;-i,-(i+1);n-1})_{m-l,m-l-1}\\
&= (-1)^{m-1}(1+w)(1+mw)u^{m+1}\tau_{-1}^{m-4}\hat{\tau}_{-2}\tilde{\tau}_{-2}\\
&\quad + (-1)^{m-2}[1+(m-1)w]u^{m-1}\tau_{-1}^{m-4}\hat{\tau}_{-2}\tilde{\tau}_{-2}\\
&\quad + \sum_{k=2}^{i-2}(-1)^{m-k-1}\tau_{m-k}|R_{\theta,n;-i,-(i+1);k}|\tau_{-1}^{m-k-3}\hat{\tau}_{-2}\tilde{\tau}_{-2}\\
&\quad + (-1)^{m-i}\hat{\tau}_{m-i}|R_{\theta,n;-i,-(i+1);i-1}|\tau_{-1}^{m-i-1}\tilde{\tau}_{-2}\\
&\quad + (-1)^{m-i-1}\tilde{\tau}_{m-i-1}|R_{\theta,n;-i,-(i+1);i}|\tau_{-1}^{m-i-1}\\
&\quad + \sum_{k=i+1}^{m-1}(-1)^{m-k-1}\tau_{m-k-1}|R_{\theta,n;-i,-(i+1);k}|\tau_{-1}^{m-k-1}\\
&= (-\tau_{-1})^m\tilde{\xi}_m\\
&\quad + (-\tau_{-1})^{m-1}\sum_{k=i+1}^{m-1}\tau_{m-k-1}(-\tau_{-1})^{-k}|R_{\theta,n;-i,-(i+1);k}|,
\end{aligned}
$$

where

$$
\begin{aligned}
\tilde{\xi}_m &= \tau_{m-1}^*\tau_{-1}^{-4}\hat{\tau}_{-2}\tilde{\tau}_{-2} + (-1)^{-i}\hat{\tau}_{m-i}|R_{\theta,n;i-1}|\tau_{-1}^{-i-1}\tilde{\tau}_{-2}\\
&\quad + (-1)^{-i-1}\tilde{\tau}_{m-i-1}|R_{\theta,n;-i,-(i+1);i}|\tau_{-1}^{-i-1}\\
&\quad + \sum_{k=2}^{i-2}(-1)^{-k-1}\tau_{m-k}|R_{\theta,n;k}|\tau_{-1}^{-k-3}\hat{\tau}_{-2}\tilde{\tau}_{-2}\\
&= \tau_{m-1}^*\tau_{-1}^{-4}\hat{\tau}_{-2}\tilde{\tau}_{-2} + (-1)^{-i}\hat{\tau}_{m-i}|R_{\theta,n;i-1}|\tau_{-1}^{-i-1}\tilde{\tau}_{-2}
\end{aligned}
$$

$$
\begin{aligned}
(20)\qquad & + (-1)^{-i-1}\tilde{\tau}_{m-i-1}|R_{\theta,n;-i,-(i+1);i}|\tau_{-1}^{-i-1} \\
& - \frac{\hat{\tau}_{-2}\tilde{\tau}_{-2}u^{m-1}(1-u^2)^2(\tilde{a}\alpha_1+u^2\tilde{b})}{\tau_{-1}^4(\alpha_1-\alpha_2)[\alpha_1(b+\tau_{-1}u)-\tau_{-1}u]}\left\{\left[\frac{\alpha_1(b+\tau_{-1}u)}{\tau_{-1}u}\right]^{i-3}-1\right\} \\
& + \frac{\hat{\tau}_{-2}\tilde{\tau}_{-2}u^{m-1}(1-u^2)^2(\tilde{a}\alpha_2+u^2\tilde{b})}{\tau_{-1}^4(\alpha_1-\alpha_2)[\alpha_2(b+\tau_{-1}u)-\tau_{-1}u]}\left\{\left[\frac{\alpha_2(b+\tau_{-1}u)}{\tau_{-1}u}\right]^{i-3}-1\right\} \\
& - \frac{\hat{\tau}_{-2}\tilde{\tau}_{-2}[(\tilde{a}\alpha_1+u^2\tilde{b})S_{i,m,1}-(\tilde{a}\alpha_2+u^2\tilde{b})S_{i,m,2}]}{\tau_{-1}^3(\alpha_1-\alpha_2)}.
\end{aligned}
$$

The last equality follows from (33). An immediate consequence is that

$$
|R_{\theta,n;-i,-(i+1);i+1}| = (-\tau_{-1})^{i+1}\tilde{\xi}_{i+1}.
$$

Now repeating the above argument, we obtain, for $4 \le i+1 < m \le n-1$,

$$
\begin{aligned}
|R_{\theta,n;-i,-(i+1);m}| &= (-\tau_{-1})^m\tilde{\xi}_m + (-\tau_{-1})^{m-1}\sum_{k_1=i+1}^{m-1}\tau_{m-k_1-1}\tilde{\xi}_{k_1} \\
&\quad + (-\tau_{-1})^{m-2}\sum_{k_1=i+2}^{m-1}\sum_{k_2=i+1}^{k_1-1}\tau_{m-k_1-1}\tau_{k_1-k_2-1}\tilde{\xi}_{k_2} \\
&\quad + \cdots + (-\tau_{-1})^{i+1}\tau_0^{m-i-1}\tilde{\xi}_{i+1}.
\end{aligned}
$$

Hence, it follows from Lemma A.4 that $|R_{\theta,n;-i,-(i+1),m}|, 4 \le i+1 \le m \le n-1$, equals the coefficient of $z^m$ in the series expansion of

$$
\frac{(-\tau_{-1}u)^m G_{\tilde{\xi}}}{(\tau_{-1}u)^{-1}zG_\tau(z)+1}.
$$

Following the proof of (19), we conclude that, for $4 \le i+1 \le m \le n-1$,

$$
\begin{aligned}
&|R_{\theta,n;-i,-(i+1),m}| \\
&\quad = \frac{(-1)^m(b+\tau_{-1}u)^m u^{2(i+1)}[a_{\tilde{\xi}}(\alpha_1^{m-i}-\alpha_2^{m-i})+u^2 b_{\tilde{\xi}}(\alpha_1^{m-i-1}-\alpha_2^{m-i-1})]}{\alpha_1-\alpha_2} \\
&\quad = \frac{(-1)^m(\tau_{-1}u)^{i+1}a_{\tilde{\xi}}}{\alpha_1-\alpha_2} \\
&\qquad \times \{(\alpha_1-u^2)[(b+\tau_{-1}u)\alpha_1]^{m-i-1}-(\alpha_2-u^2)[(b+\tau_{-1}u)\alpha_2]^{m-i-1}\} \\
&\qquad + \frac{(-1)^m(\tau_{-1}u)^{i+1}u^2(a_{\tilde{\xi}}+b_{\tilde{\xi}})}{\alpha_1-\alpha_2} \\
&\qquad\quad \times \{[(b+\tau_{-1}u)\alpha_1]^{m-i-1}-[(b+\tau_{-1}u)\alpha_2]^{m-i-1}\},
\end{aligned}
$$

and (15) follows from Lemma A.4 after some algebra. Finally, we observe from Case 3 that, for $5 \le i+2 \le j \le n$,

$$|R_{\theta,n;-i,-j;n-1}| = \tau_{-1}^{j-i-2} |R_{\theta,n;-i,-j;i+1}| |\tilde{A}_{\theta,n;n-j}|,$$

and (16) follows from Lemmas A.2 and A.3 (see Appendix A). This proves Theorem 4. □

Remark. We wish to add that the exact results of Sections 2 and 3 have been checked by *Mathematica* [15].

**4. Asymptotic approximations.** This section establishes asymptotic approximations for $|R_{\theta,n}|$ and $R_{\theta,n}^{-1}$. These approximations are very sharp in that the error is of the order $O((2+\sqrt{3})^{-n})$ as $n \to \infty$. The following theorem gives the approximation for the determinant.

Theorem 5. *Let $0 < \beta_0 < \beta_1 < \infty$. Then*

$$\begin{aligned}|R_{\theta,n}| &= \frac{\tilde{a}\alpha_1 + u^2\tilde{b}}{\alpha_1 - \alpha_2}[-(b+\tau_{-1}u)\alpha_1]^{n-2}\left[1 + O\left(\left(\frac{2-\sqrt{3}}{2+\sqrt{3}}\right)^n\right)\right] \\ &= \frac{w^{3n-4}e^{-2(n-2)w}\sqrt{3}}{4}\left[\frac{2(2+\sqrt{3})}{3}\right]^{n-1}\left[1 + \frac{(12+7\sqrt{3})nw^2}{60+30\sqrt{3}} + O(w^2)\right],\end{aligned}$$

*as $n \to \infty$ uniformly over $\theta \in [\beta_0, \beta_1]$.*

Proof. From the definitions of $\alpha_1$ and $\alpha_2$ in Theorem 1 and using the *Mathematica* command

$$\texttt{Simplify}[\texttt{Series}[\alpha_1/\alpha_2, \{\texttt{w}, 0, 0\}]],$$

we observe that

$$\frac{\alpha_2}{\alpha_1} = \frac{2-\sqrt{3}}{2+\sqrt{3}}[1 + O(w^2)], \tag{21}$$

as $n \to \infty$ uniformly over $\theta \in [\beta_0, \beta_1]$. Using *Mathematica* in a similar way, we also have

$$\frac{\tilde{a}\alpha_2 + u^2\tilde{b}}{\tilde{a}\alpha_1 + u^2\tilde{b}} = \frac{2-\sqrt{3}}{2+\sqrt{3}}[1 + O(w)],$$

as $n \to \infty$ uniformly over $\theta \in [\beta_0, \beta_1]$. Thus, it follows from Theorem 1 that

$$|R_{\theta,n}| = \frac{\tilde{a}\alpha_1 + u^2\tilde{b}}{\alpha_1 - \alpha_2}[-(b+\tau_{-1}u)\alpha_1]^{n-2}\left[1 + O\left(\left(\frac{2-\sqrt{3}}{2+\sqrt{3}}\right)^n\right)\right],$$

as $n \to \infty$ uniformly over $\theta \in [\beta_0, \beta_1]$. Next, using the notation of Theorem 1, we observe that

$$-e^{4w}(a - 2\tau_{-1}u) = 1 + 4we^{2w} - e^{4w} = -\frac{8w^3}{3}\sum_{k=0}^{\infty}\frac{3!2^k(2^{k+2}-k-3)w^k}{(k+3)!}$$

and

$$e^{8w}[(a-2\tau_{-1}u)^2 - 4\tau_{-1}u(b+\tau_{-1}u)] = \frac{16w^6}{3}(1-\zeta_w),$$

where

$$\begin{aligned}\zeta_w &= 96w\sum_{k=0}^{\infty}\frac{(2w)^k}{(k+7)!} - 18432w\sum_{k=0}^{\infty}\frac{(4w)^k}{(k+7)!} + 209952w\sum_{k=0}^{\infty}\frac{(6w)^k}{(k+7)!}\\ &\quad - 393216w\sum_{k=0}^{\infty}\frac{(8w)^k}{(k+7)!} + 24w\sum_{k=0}^{\infty}\frac{(2w)^k}{(k+5)!}\\ &\quad - 1536w\sum_{k=0}^{\infty}\frac{(4w)^k}{(k+5)!} + 5832w\sum_{k=0}^{\infty}\frac{(6w)^k}{(k+5)!}\\ &= -4w(1+O(w)),\end{aligned}$$

as $n \to \infty$ uniformly over $\theta \in [\beta_0, \beta_1]$. Hence,

$$\begin{aligned}&(b+\tau_{-1}u)\alpha_1\\ &\quad = -\frac{e^{-4w}}{2}[e^{4w}(a-2\tau_{-1}u) + e^{4w}\sqrt{(a-2\tau_{-1}u)^2 - 4\tau_{-1}u(b+\tau_{-1}u)}\,]\\ &\quad = -e^{-4w}\left[\frac{4w^3}{3}\sum_{k=0}^{\infty}\frac{3!2^k(2^{k+2}-k-3)w^k}{(k+3)!} + \frac{2w^3}{\sqrt{3}}(1-\zeta_w)^{1/2}\right] \qquad (22)\\ &\quad = -e^{-4w}\left\{\frac{4w^3}{3}\sum_{k=0}^{\infty}\frac{3!2^k(2^{k+2}-k-3)w^k}{(k+3)!}\right.\\ &\qquad\qquad \left. + \frac{2w^3}{\sqrt{3}}\left[1 - \frac{\zeta_w}{2} - \zeta_w^2\sum_{k=0}^{\infty}\frac{(2k+1)!\zeta_w^k}{4^{k+1}k!(k+2)!}\right]\right\}\\ &\quad = -\frac{2(2+\sqrt{3})w^3}{3}(1+\Delta_{1,w}),\end{aligned}$$

where

$$\begin{aligned}\Delta_{1,w} &= \frac{2we^{-4w}}{2+\sqrt{3}}\sum_{k=0}^{\infty}\frac{3!2^{k+1}(2^{k+3}-k-4)w^k}{(k+4)!} - 4we^{-4w}\sum_{k=0}^{\infty}\frac{(4w)^k}{(k+1)!}\\ &\quad - \frac{\sqrt{3}\zeta_w e^{-4w}}{2+\sqrt{3}}\left[\frac{1}{2} + \zeta_w\sum_{k=0}^{\infty}\frac{(2k+1)!\zeta_w^k}{4^{k+1}k!(k+2)!}\right]\end{aligned}$$

$$= -2w + \frac{(132+67\sqrt{3})w^2}{60+30\sqrt{3}} + O(w^3),$$

as $n \to \infty$ uniformly over $\theta \in [\beta_0, \beta_1]$. We conclude from (22) and Theorem 1 that

$$\begin{aligned}
|R_{\theta,n}| &= \frac{(2+\sqrt{3})w^{3n-4}}{2\sqrt{3}}\left[\frac{2(2+\sqrt{3})(1+\Delta_{1,w})}{3}\right]^{n-2} \\
&\quad \times \left[1 - \frac{2(3+\sqrt{3})w}{3(2+\sqrt{3})} + O(w^2)\right] \\
&= \frac{w^{3n-4}e^{-2(n-2)w}\sqrt{3}}{4}\left[\frac{2(2+\sqrt{3})}{3}\right]^{n-1} \\
&\quad \times \left[1 + \frac{(12+7\sqrt{3})nw^2}{60+30\sqrt{3}} + O(w^2)\right],
\end{aligned}$$

as $n \to \infty$ uniformly over $\theta \in [\beta_0, \beta_1]$. This proves Theorem 5. □

The next two theorems give the approximation for the inverse $R_{\theta,n}^{-1}$.

THEOREM 6. *Let $0 < \beta_0 < \beta_1 < \infty$. Then with the notation of Appendix B,*

$$\begin{aligned}
(R_{\theta,n}^{-1})_{2,n-1} &= (-1)^{n-3}\frac{\tau_0\{\tilde{\tau}_{-2}[(1+2w)^2u^4-1] + \tilde{\tau}_{-1}\tau_1^*\}(\alpha_1-\alpha_2)}{\tau_{-1}^3(\tilde{a}\alpha_1 + u^2\tilde{b})} \\
&\quad \times \left(-\frac{u}{\alpha_1}\right)^{n-2}\left[1 + O\left(\left(\frac{2-\sqrt{3}}{2+\sqrt{3}}\right)^n\right)\right],
\end{aligned}$$

$$\begin{aligned}
(R_{\theta,n}^{-1})_{n-1,n-1} &= \left\{\frac{C_{7,1}}{[-(b+\tau_{-1}u)\alpha_1]^2} - \frac{C_{7,2}(\tilde{a}\alpha_2 + u^2\tilde{b})}{(\tilde{a}\alpha_1 + u^2\tilde{b})[-(b+\tau_{-1}u)\alpha_1]^2}\left(\frac{\alpha_2}{\alpha_1}\right)^{n-4}\right\} \\
&\quad \times \left[1 + O\left(\left(\frac{2-\sqrt{3}}{2+\sqrt{3}}\right)^n\right)\right],
\end{aligned}$$

$$(R_{\theta,n}^{-1})_{1,n} = (-1)^{n+1}\frac{w^2u^2(\alpha_1-\alpha_2)}{\tau_{-1}(\tilde{a}\alpha_1 + u^2\tilde{b})}\left(-\frac{u}{\alpha_1}\right)^{n-2}\left[1 + O\left(\left(\frac{2-\sqrt{3}}{2+\sqrt{3}}\right)^n\right)\right],$$

$$\begin{aligned}
(R_{\theta,n}^{-1})_{2,n} &= (-1)^{n+2}\frac{\{\tilde{\tau}_{-2}[(1+2w)^2u^4-1] + \tilde{\tau}_{-1}\tau_1^*\}(\alpha_1-\alpha_2)}{\tau_{-1}^2(\tilde{a}\alpha_1 + u^2\tilde{b})}\left(-\frac{u}{\alpha_1}\right)^{n-2} \\
&\quad \times \left[1 + O\left(\left(\frac{2-\sqrt{3}}{2+\sqrt{3}}\right)^n\right)\right],
\end{aligned}$$

$$(R_{\theta,n}^{-1})_{i,n} = \frac{(-1)^{n+i}C_{6,1}\tau_{-1}^{n-i-2}[-(b+\tau_{-1}u)\alpha_1]^{i-3}}{[-(b+\tau_{-1}u)\alpha_1]^{n-2}}\left[1 + O\left(\left(\frac{2-\sqrt{3}}{2+\sqrt{3}}\right)^n\right)\right]$$

$$
\begin{aligned}
&\quad - \frac{(-1)^{n+i} C_{6,2} \tau_{-1}^{n-i-2} (\tilde{a}\alpha_2 + u^2 \tilde{b})[-(b+\tau_{-1}u)\alpha_2]^{i-3}}{(\tilde{a}\alpha_1 + u^2\tilde{b})[-(b+\tau_{-1}u)\alpha_1]^{n-2}} \\
&\quad \times \left[1 + O\left(\left(\frac{2-\sqrt{3}}{2+\sqrt{3}}\right)^n\right)\right], \\
(R_{\theta,n}^{-1})_{n-1,n} &= -\frac{C_{3,1}}{[-(b+\tau_{-1}u)\alpha_1]^2}\left[1 + O\left(\left(\frac{2-\sqrt{3}}{2+\sqrt{3}}\right)^n\right)\right], \\
(R_{\theta,n}^{-1})_{n,n} &= -\frac{1}{(b+\tau_{-1}u)\alpha_1}\left[1 + O\left(\left(\frac{2-\sqrt{3}}{2+\sqrt{3}}\right)^n\right)\right], \\
(R_{\theta,n}^{-1})_{1,1} &= -\frac{1}{(b+\tau_{-1}u)\alpha_1}\left[1 + O\left(\left(\frac{2-\sqrt{3}}{2+\sqrt{3}}\right)^n\right)\right], \\
(R_{\theta,n}^{-1})_{2,1} &= -\frac{C_{3,1}}{[-(b+\tau_{-1}u)\alpha_1]^2}\left[1 + O\left(\left(\frac{2-\sqrt{3}}{2+\sqrt{3}}\right)^n\right)\right], \\
(R_{\theta,n}^{-1})_{i,1} &= \frac{(-1)^{i+1} C_{6,1}\tau_{-1}^{i-3}[-(b+\tau_{-1}u)\alpha_1]^{n-i-2}}{[-(b+\tau_{-1}u)\alpha_1]^{n-2}}\left[1 + O\left(\left(\frac{2-\sqrt{3}}{2+\sqrt{3}}\right)^n\right)\right] \\
&\quad - \frac{(-1)^{i+1} C_{6,2}\tau_{-1}^{i-3}(\tilde{a}\alpha_2 + u^2\tilde{b})[-(b+\tau_{-1}u)\alpha_2]^{n-i-2}}{(\tilde{a}\alpha_1 + u^2\tilde{b})[-(b+\tau_{-1}u)\alpha_1]^{n-2}} \\
&\quad \times \left[1 + O\left(\left(\frac{2-\sqrt{3}}{2+\sqrt{3}}\right)^n\right)\right]
\end{aligned}
$$

*and*

$$
\begin{aligned}
(R_{\theta,n}^{-1})_{n-1,1} &= (-1)^{n+2}\frac{\{\tilde{\tau}_{-2}[(1+2w)^2u^4 - 1] + \tilde{\tau}_{-1}\tau_1^*\}(\alpha_1 - \alpha_2)}{\tau_{-1}^2(\tilde{a}\alpha_1 + u^2\tilde{b})}\left(-\frac{u}{\alpha_1}\right)^{n-2} \\
&\quad \times \left[1 + O\left(\left(\frac{2-\sqrt{3}}{2+\sqrt{3}}\right)^n\right)\right],
\end{aligned}
$$

*as $n \to \infty$ uniformly over $\theta \in [\beta_0, \beta_1]$ and $3 \le i \le n-2$.*

PROOF. Theorem 6 is a consequence of (10), (21) and Theorems 1–4. □

THEOREM 7. *Let $0 < \beta_0 < \beta_1 < \infty$. Then with the notation of Appendix B,*

$$
\begin{aligned}
(R_{\theta,n}^{-1})_{i,i} &= \left\{\frac{C_{1,1}\tau_{-1}u^4}{\alpha_1 - \alpha_2}\left\{\frac{\alpha_1 - u^2}{[-(b+\tau_{-1}u)\alpha_1]^3} - \frac{\alpha_2 - u^2}{[-(b+\tau_{-1}u)\alpha_1]^3}\left(\frac{\alpha_2}{\alpha_1}\right)^{n-i-2}\right\}\right. \\
&\quad - \frac{C_{1,2}\tau_{-1}u^4(\tilde{a}\alpha_2 + u^2\tilde{b})}{(\alpha_1 - \alpha_2)(\tilde{a}\alpha_1 + u^2\tilde{b})}
\end{aligned}
$$

$$\times\left\{\frac{\alpha_1-u^2}{[-(b+\tau_{-1}u)\alpha_2]^3}\left(\frac{\alpha_2}{\alpha_1}\right)^{i}-\frac{\alpha_2-u^2}{[-(b+\tau_{-1}u)\alpha_2]^3}\left(\frac{\alpha_2}{\alpha_1}\right)^{n-2}\right\}$$

$$+\frac{C_{2,1}\tau_{-1}u^6}{\alpha_1-\alpha_2}\left\{\frac{1}{[-(b+\tau_{-1}u)\alpha_1]^3}-\frac{1}{[-(b+\tau_{-1}u)\alpha_1]^3}\left(\frac{\alpha_2}{\alpha_1}\right)^{n-i-2}\right\}$$

$$-\frac{C_{2,2}\tau_{-1}u^6(\tilde{a}\alpha_2+u^2\tilde{b})}{(\alpha_1-\alpha_2)(\tilde{a}\alpha_1+u^2\tilde{b})}$$

$$\times\left\{\frac{1}{[-(b+\tau_{-1}u)\alpha_2]^3}\left(\frac{\alpha_2}{\alpha_1}\right)^{i}-\frac{1}{[-(b+\tau_{-1}u)\alpha_2]^3}\left(\frac{\alpha_2}{\alpha_1}\right)^{n-2}\right\}\Bigg\}$$

$$\times\left[1+O\left(\left(\frac{2-\sqrt{3}}{2+\sqrt{3}}\right)^{n}\right)\right],$$

*as* $n\to\infty$ *uniformly over* $\theta\in[\beta_0,\beta_1]$ *and* $3\le i\le n-2$,

$$(R_{\theta,n}^{-1})_{i,j}=\frac{1}{\tau_{-1}^4(\alpha_1-\alpha_2)}\left(\frac{u}{\alpha_1}\right)^{j-i+2}$$

$$\times\left\{C_{6,1}\left[(a\alpha_1+u^2b)-(a\alpha_2+u^2b)\left(\frac{\alpha_2}{\alpha_1}\right)^{n-j-1}\right]\right.$$

$$-\frac{C_{6,2}(\tilde{a}\alpha_2+u^2\tilde{b})}{\tilde{a}\alpha_1+u^2\tilde{b}}$$

$$\left.\times\left[(a\alpha_1+u^2b)\left(\frac{\alpha_2}{\alpha_1}\right)^{i-3}-(a\alpha_2+u^2b)\left(\frac{\alpha_2}{\alpha_1}\right)^{n-j+i-4}\right]\right\}$$

$$\times\left[1+O\left(\left(\frac{2-\sqrt{3}}{2+\sqrt{3}}\right)^{n}\right)\right],$$

*as* $n\to\infty$ *uniformly over* $\theta\in[\beta_0,\beta_1]$ *and* $5\le i+2\le j\le n-1$, *and*

$$(R_{\theta,n}^{-1})_{i,i+1}$$

$$=-\left\{\frac{C_{4,1}\tau_{-1}u^4}{\alpha_1-\alpha_2}\right.$$

$$\times\left\{\frac{\alpha_1-u^2}{[-(b+\tau_{-1}u)\alpha_1]^3}-\frac{\alpha_2-u^2}{[-(b+\tau_{-1}u)\alpha_1]^3}\left(\frac{\alpha_2}{\alpha_1}\right)^{n-i-2}\right\}$$

$$-\frac{C_{4,2}\tau_{-1}u^4(\tilde{a}\alpha_2+u^2\tilde{b})}{(\alpha_1-\alpha_2)(\tilde{a}\alpha_1+u^2\tilde{b})}$$

$$\times\left\{\frac{\alpha_1-u^2}{[-(b+\tau_{-1}u)\alpha_2]^3}\left(\frac{\alpha_2}{\alpha_1}\right)^{i}-\frac{\alpha_2-u^2}{[-(b+\tau_{-1}u)\alpha_2]^3}\left(\frac{\alpha_2}{\alpha_1}\right)^{n-2}\right\}$$

$$
\begin{aligned}
&+ \frac{C_{5,1}\tau_{-1}u^6}{\alpha_1-\alpha_2} \\
&\quad\times \biggl\{\frac{1}{[-(b+\tau_{-1}u)\alpha_1]^3} - \frac{1}{[-(b+\tau_{-1}u)\alpha_1]^3}\biggl(\frac{\alpha_2}{\alpha_1}\biggr)^{n-i-2}\biggr\} \\
&- \frac{C_{5,2}\tau_{-1}u^6(\tilde a\alpha_2+u^2\tilde b)}{(\alpha_1-\alpha_2)(\tilde a\alpha_1+u^2\tilde b)} \\
&\quad\times \biggl\{\frac{1}{[-(b+\tau_{-1}u)\alpha_2]^3}\biggl(\frac{\alpha_2}{\alpha_1}\biggr)^{i} - \frac{1}{[-(b+\tau_{-1}u)\alpha_2]^3}\biggl(\frac{\alpha_2}{\alpha_1}\biggr)^{n-2}\biggr\}\biggr\} \\
\times &\biggl[1+O\biggl(\biggl(\frac{2-\sqrt{3}}{2+\sqrt{3}}\biggr)^n\biggr)\biggr],
\end{aligned}
$$

*as* $n\to\infty$ *uniformly over* $\theta\in[\beta_0,\beta_1]$ *and* $3\le i\le n-2$. *Finally,*

$$
\begin{aligned}
(R_{\theta,n}^{-1})_{i,n} = &\biggl\{\frac{C_{6,1}}{[-(b+\tau_{-1}u)\alpha_1]^3}\biggl(\frac{u}{\alpha_1}\biggr)^{n-i-2} \\
&- \frac{C_{6,2}(\tilde a\alpha_2+u^2\tilde b)}{(\tilde a\alpha_1+u^2\tilde b)[-(b+\tau_{-1}u)\alpha_2]^3}\biggl(\frac{u}{\alpha_1}\biggr)^{n-i-2}\biggl(\frac{\alpha_2}{\alpha_1}\biggr)^{i}\biggr\} \\
&\times\biggl[1+O\biggl(\biggl(\frac{2-\sqrt{3}}{2+\sqrt{3}}\biggr)^n\biggr)\biggr],
\end{aligned}
$$

*as* $n\to\infty$ *uniformly over* $\theta\in[\beta_0,\beta_1]$ *and* $3\le i\le n-2$.

PROOF. As in the proof of Theorem 6, Theorem 7 is a consequence of (10), (21) and Theorems 1–4. □

**5. Fisher information.** Stein [12], page 178, noted that the asymptotic theory of maximum likelihood estimation based on independently and identically distributed observations suggests that calculating the Fisher information matrix is a fruitful way of learning about the behavior of the maximum likelihood estimators. Abt and Welch [1] showed via three examples (either proven analytically or justified by simulation) that, for the covariance parameters of Gaussian processes under fixed-domain asymptotics, the covariance matrix of the limiting distribution of their maximum likelihood estimators equals the limit of the inverse Fisher information matrix. Also, Section 6.6 of [12] reported a numerical study of the Fisher information matrix for observations from a mean 0 Gaussian process on the real line with a Matérn-type covariance function.

Let $\tilde X_n$ be the random vector as in Section 1, having a $n^d$-variate normal distribution with mean 0 and covariance matrix $\Sigma_{\phi,\theta_1,\dots,\theta_d;n}$. Motivated by the work described in the previous paragraph, the following theorem computes the Fisher information matrix for the parameters $\phi, \theta_1, \dots, \theta_d$ based on $\tilde{X}_n$.

THEOREM 8. *Let $\tilde{X}_n \sim N_{n^d}(0, \Sigma_{\phi,\theta_1,\dots,\theta_d})$. Then the elements of the Fisher information matrix for the parameters $\phi, \theta_1, \dots, \theta_d$ are given by*

$$-E\left[\frac{\partial^2}{\partial\phi^2}\log L_n(\phi,\theta_1,\dots,\theta_d)\right] = \frac{d^2 n^d}{2\phi^2},$$

$$-E\left[\frac{\partial^2}{\partial\theta_t^2}\log L_n(\phi,\theta_1,\dots,\theta_d)\right] = \frac{n^{d-1}(2\theta_t+5)}{\theta_t^2} + O(n^{d-2}),$$

$$-E\left[\frac{\partial^2}{\partial\phi\,\partial\theta_t}\log L_n(\phi,\theta_1,\dots,\theta_d)\right] = -\frac{dn^{d-1}(\theta_t+2)}{\phi\theta_t} + O(n^{d-2}),$$

$$-E\left[\frac{\partial^2}{\partial\theta_s\,\partial\theta_t}\log L_n(\phi,\theta_1,\dots,\theta_d)\right] = O(n^{d-2}),$$

*as $n \to \infty$ whenever $1 \le s \ne t \le d$.*

PROOF. We observe from [2], page 600, that

$$\left|\bigotimes_{t=1}^{d} R_{\theta_t,n}\right| = \prod_{t=1}^{d}[|R_{\theta_t,n}|^{n^{d-1}}],$$

and, hence, it follows from (5) that

$$\begin{aligned} &2\log L_n(\phi,\theta_1,\dots,\theta_d) \\ &\quad= -n^d\log(2\pi) - dn^d\log\left(\frac{\pi}{2}\right) - dn^d\log(\phi) + 3n^d\sum_{t=1}^{d}\log(\theta_t) \\ &\qquad - n^{d-1}\sum_{t=1}^{d}\log(|R_{\theta_t,n}|) - \frac{2^d\theta_1^3\cdots\theta_d^3}{\pi^d\phi^d}\tilde{X}_n'\left(\bigotimes_{t=1}^{d}R_{\theta_t,n}\right)^{-1}\tilde{X}_n. \end{aligned} \tag{23}$$

Consequently,

$$-E\left[\frac{\partial^2}{\partial\phi^2}\log L_n(\phi,\theta_1,\dots,\theta_d)\right] = -\frac{dn^d}{2\phi^2} + \frac{d(d+1)n^d}{2\phi^2} = \frac{d^2n^d}{2\phi^2}.$$

Writing $w_t = \theta_t/n$, we observe that

$$\begin{aligned} &\frac{\partial^2}{\partial\theta_t^2}[3n^d\log(\theta_t)] = -\frac{3n^d}{\theta_t^2}, \\ &\frac{\partial^2}{\partial\theta_t^2}(n^{d-1}\log|R_{\theta_t,n}|) = \frac{\partial}{\partial\theta_t}\left(n^{d-1}|R_{\theta_t,n}|^{-1}\frac{\partial|R_{\theta_t,n}|}{\partial\theta_t}\right) \\ &\qquad = \frac{n^{d-3}}{|R_{\theta_t,n}|}\frac{\partial^2|R_{\theta_t,n}|}{\partial w_t^2} - \frac{n^{d-3}}{|R_{\theta_t,n}|^2}\left(\frac{\partial|R_{\theta_t,n}|}{\partial w_t}\right)^2. \end{aligned} \tag{24}$$

Next, we observe from Lemma A.6 that, as $n \to \infty$,

$$\begin{aligned} \frac{1}{n}\operatorname{tr}\Big[\Big(\frac{\partial}{\partial\theta_d}R^{-1}_{\theta_d,n}\Big)R_{\theta_d,n}\Big] &= -\frac{3}{\theta_d} + \frac{2(\theta_d+2)}{\theta_d n} + O(n^{-2}), \\ \frac{1}{n}\operatorname{tr}\Big[\Big(\frac{\partial^2}{\partial\theta_d^2}R^{-1}_{\theta_d,n}\Big)R_{\theta_d,n}\Big] &= \frac{12}{\theta_d^2} - \frac{2(4\theta_d+9)}{\theta_d^2 n} + O(n^{-2}). \end{aligned} \tag{25}$$

We observe that

$$\begin{aligned} &\frac{\partial^2}{\partial\theta_d^2}\Bigg[\theta_1^3\cdots\theta_d^3\tilde{X}_n'\Bigg(\bigotimes_{t=1}^{d}R_{\theta_t,n}\Bigg)^{-1}\tilde{X}_n\Bigg] \\ &\quad = \frac{\partial}{\partial\theta_d}\Bigg\{3\theta_1^3\cdots\theta_d^2\tilde{X}_n'\Bigg(\bigotimes_{t=1}^{d}R^{-1}_{\theta_t,n}\Bigg)\tilde{X}_n \\ &\qquad\qquad + \theta_1^3\cdots\theta_d^3\tilde{X}_n'\Bigg[\Bigg(\bigotimes_{t=1}^{d-1}R^{-1}_{\theta_t,n}\Bigg)\otimes\frac{\partial}{\partial\theta_d}(R^{-1}_{\theta_d,n})\Bigg]\tilde{X}_n\Bigg\} \\ &\quad = 6\theta_1^3\cdots\theta_d\tilde{X}_n'\Bigg(\bigotimes_{t=1}^{d}R^{-1}_{\theta_t,n}\Bigg)\tilde{X}_n \\ &\qquad + 6\theta_1^3\cdots\theta_d^2\tilde{X}_n'\Bigg[\Bigg(\bigotimes_{t=1}^{d-1}R^{-1}_{\theta_t,n}\Bigg)\otimes\frac{\partial}{\partial\theta_d}(R^{-1}_{\theta_d,n})\Bigg]\tilde{X}_n \\ &\qquad + \theta_1^3\cdots\theta_d^3\tilde{X}_n'\Bigg[\Bigg(\bigotimes_{t=1}^{d-1}R^{-1}_{\theta_t,n}\Bigg)\otimes\frac{\partial^2}{\partial\theta_d^2}(R^{-1}_{\theta_d,n})\Bigg]\tilde{X}_n, \end{aligned}$$

and, hence, using (25),

$$\begin{aligned} &E\frac{\partial^2}{\partial\theta_d^2}\Bigg[\frac{2^d\theta_1^3\cdots\theta_d^3}{\pi^d\phi^d}\tilde{X}_n'\Bigg(\bigotimes_{t=1}^{d}R_{\theta_t,n}\Bigg)^{-1}\tilde{X}_n\Bigg] \\ &\quad = \frac{6(2^d)\theta_1^3\cdots\theta_d}{\pi^d\phi^d}\operatorname{tr}\Bigg[\Bigg(\bigotimes_{t=1}^{d}R^{-1}_{\theta_t,n}\Bigg)\Sigma_{\phi,\theta_1,\ldots,\theta_d;n}\Bigg] \\ &\qquad + \frac{6(2^d)\theta_1^3\cdots\theta_d^2}{\pi^d\phi^d}\operatorname{tr}\Bigg\{\Bigg[\Bigg(\bigotimes_{t=1}^{d-1}R^{-1}_{\theta_t,n}\Bigg)\otimes\frac{\partial}{\partial\theta_d}(R^{-1}_{\theta_d,n})\Bigg]\Sigma_{\phi,\theta_1,\ldots,\theta_d;n}\Bigg\} \\ &\qquad + \frac{2^d\theta_1^3\cdots\theta_d^3}{\pi^d\phi^d}\operatorname{tr}\Bigg\{\Bigg[\Bigg(\bigotimes_{t=1}^{d-1}R^{-1}_{\theta_t,n}\Bigg)\otimes\frac{\partial^2}{\partial\theta_d^2}(R^{-1}_{\theta_d,n})\Bigg]\Sigma_{\phi,\theta_1,\ldots,\theta_d;n}\Bigg\} \\ &\quad = \frac{6n^d}{\theta_d^2} + \frac{6n^{d-1}}{\theta_d}\operatorname{tr}\Big[\Big(\frac{\partial}{\partial\theta_d}R^{-1}_{\theta_d,n}\Big)R_{\theta_d,n}\Big] + n^{d-1}\operatorname{tr}\Big[\Big(\frac{\partial^2}{\partial\theta_d^2}R^{-1}_{\theta_d,n}\Big)R_{\theta_d,n}\Big] \\ &\quad = \frac{2n^{d-1}(2\theta_d+3)}{\theta_d^2} + O(n^{d-2}). \end{aligned} \tag{26}$$

We conclude from (23), (24), (26) and Theorem 5 that, as $n\to\infty$,

$$
\begin{aligned}
-E&\frac{\partial^2}{\partial\theta_d^2}\log L_n(\phi,\theta_1,\ldots,\theta_d)\\
&=\frac{3n^d}{2\theta_d^2}-\left(\frac{3n^d}{2\theta_d^2}-\frac{2n^{d-1}}{\theta_d^2}\right)+\frac{n^{d-1}(2\theta_d+3)}{\theta_d^2}+O(n^{d-2})\\
&=\frac{n^{d-1}(2\theta_d+5)}{\theta_d^2}+O(n^{d-2}).
\end{aligned}
$$

Finally, using (25) again, we have, as $n\to\infty$,

$$
\begin{aligned}
-E&\frac{\partial^2}{\partial\theta_d\,\partial\phi}\log L_n(\phi,\theta_1,\ldots,\theta_d)\\
&=-E\frac{\partial}{\partial\theta_d}\left[\frac{2^{d-1}d\theta_1^3\cdots\theta_d^3}{\pi^d\phi^{d+1}}\tilde X_n'\left(\bigotimes_{t=1}^d R_{\theta_t,n}\right)^{-1}\tilde X_n\right]\\
&=-\frac{3(2^{d-1})d\theta_1^3\cdots\theta_d^2}{\pi^d\phi^{d+1}}\operatorname{tr}\left[\left(\bigotimes_{t=1}^d R_{\theta_t,n}\right)^{-1}\Sigma_{\phi,\theta_1,\ldots,\theta_d;n}\right]\\
&\qquad-\frac{2^{d-1}d\theta_1^3\cdots\theta_d^3}{\pi^d\phi^{d+1}}\operatorname{tr}\left\{\left[\left(\bigotimes_{t=1}^{d-1} R_{\theta_t,n}^{-1}\right)\otimes\left(\frac{\partial}{\partial\theta_d}R_{\theta_d,n}^{-1}\right)\right]\Sigma_{\phi,\theta_1,\ldots,\theta_d;n}\right\}\\
&=-\frac{3dn^d}{2\phi\theta_d}-\frac{dn^{d-1}}{2\phi}\operatorname{tr}\left[\left(\frac{\partial}{\partial\theta_d}R_{\theta_d,n}^{-1}\right)R_{\theta_d,n}\right]\\
&=-\frac{dn^{d-1}(\theta_d+2)}{\phi\theta_d}+O(n^{d-2})
\end{aligned}
$$

and

$$
\begin{aligned}
-E&\frac{\partial^2}{\partial\theta_{d-1}\,\partial\theta_d}\log L_n(\phi,\theta_1,\ldots,\theta_d)\\
&=E\frac{\partial}{\partial\theta_{d-1}}\left[\frac{3(2^{d-1})\theta_1^3\cdots\theta_d^2}{\pi^d\phi^d}\tilde X_n'\left(\bigotimes_{t=1}^d R_{\theta_t,n}^{-1}\right)\tilde X_n\right]\\
&\qquad+E\frac{\partial}{\partial\theta_{d-1}}\left[\frac{2^{d-1}\theta_1^3\cdots\theta_d^3}{\pi^d\phi^d}\tilde X_n'\left(\bigotimes_{t=1}^{d-1} R_{\theta_t,n}^{-1}\right)\otimes\left(\frac{\partial}{\partial\theta_d}R_{\theta_d,n}^{-1}\right)\tilde X_n\right]\\
&=\frac{9n^d}{2\theta_{d-1}\theta_d}+\frac{3n^{d-1}}{2\theta_d}\operatorname{tr}\left[\left(\frac{\partial}{\partial\theta_{d-1}}R_{\theta_{d-1},n}^{-1}\right)R_{\theta_{d-1},n}\right]\\
&\qquad+\frac{3n^{d-1}}{2\theta_{d-1}}\operatorname{tr}\left[\left(\frac{\partial}{\partial\theta_d}R_{\theta_d,n}^{-1}\right)R_{\theta_d,n}\right]
\end{aligned}
$$

$$+\frac{n^{d-2}}{2}\operatorname{tr}\left[\left(\frac{\partial}{\partial\theta_{d-1}}R^{-1}_{\theta_{d-1},n}\right)R_{\theta_{d-1},n}\right]\operatorname{tr}\left[\left(\frac{\partial}{\partial\theta_d}R^{-1}_{\theta_d,n}\right)R_{\theta_d,n}\right]$$
$$=O(n^{d-2}).$$

The proof of Theorem 8 is completed by invoking the symmetry of $\theta_1,\ldots,\theta_d$. $\square$

**6. Estimating the scale parameter.** Let $\tilde{X}_n \sim N_{n^d}(0,\Sigma_{\phi,\theta_1,\ldots,\theta_d;n})$ be as in Section 1. Also let $0<\beta_{0,t}<\beta_{1,t}<\infty, t=0,\ldots,d$, be known constants such that $\beta_{0,0}<\phi<\beta_{1,0}$ and $\beta_{0,t}<\theta_t<\beta_{1,t}, t=1,\ldots,d$. This section will be concerned with the construction of a consistent maximum likelihood-type estimator for the scale parameter $\phi$. First set the estimator for $(\theta_1,\ldots,\theta_d)$ to be a known arbitrary but fixed vector, say $(\tilde{\theta}_1,\ldots,\tilde{\theta}_d)\in\prod_{t=1}^d(\beta_{0,t},\beta_{1,t})$, and define the estimator $\hat{\phi}_{\tilde{\theta}_1,\ldots,\tilde{\theta}_d}$ for $\phi$ to be that value of $\phi$ that maximizes the log-likelihood function $\log L_n(\phi,\tilde{\theta}_1,\ldots,\tilde{\theta}_d)$ [see (23)]. On differentiating the log-likelihood with respect to $\phi$, we obtain

$$\frac{d}{d\phi}[2\log L_n(\phi,\tilde{\theta}_1,\ldots,\tilde{\theta}_d)]=-\frac{dn^d}{\phi}+\frac{d2^d\tilde{\theta}_1^3\cdots\tilde{\theta}_d^3}{\pi^d\phi^{d+1}}\tilde{X}_n'\left(\bigotimes_{t=1}^d R_{\tilde{\theta}_t,n}\right)^{-1}\tilde{X}_n.$$

Equating the right-hand side of the above equation to zero and then solving for $\phi$ gives us a maximum likelihood-type estimator for $\phi$, namely,

$$\hat{\phi}_{\tilde{\theta}_1,\ldots,\tilde{\theta}_d}=\left\{\frac{2^d\tilde{\theta}_1^3\cdots\tilde{\theta}_d^3}{\pi^d n^d}\tilde{X}_n'\left(\bigotimes_{t=1}^d R_{\tilde{\theta}_t,n}\right)^{-1}\tilde{X}_n\right\}^{1/d}.$$

THEOREM 9. *With the above notation and conditions, we have*

$$\frac{2^d\tilde{\theta}_1^3\cdots\tilde{\theta}_d^3}{\pi^d n^d}E\left[\tilde{X}_n'\left(\bigotimes_{t=1}^d R_{\tilde{\theta}_t,n}\right)^{-1}\tilde{X}_n\right]=\phi^d+O(n^{-1}),$$
$$\operatorname{Var}\left\{\frac{2^d\tilde{\theta}_1^3\cdots\tilde{\theta}_d^3}{\pi^d n^d}\tilde{X}_n'\left(\bigotimes_{t=1}^d R_{\tilde{\theta}_t,n}\right)^{-1}\tilde{X}_n\right\}=O(n^{-d}),$$

*as $n\to\infty$. Consequently, $\hat{\phi}_{\tilde{\theta}_1,\ldots,\tilde{\theta}_d}\to\phi$ in probability as $n\to\infty$.*

PROOF. Since $(\tilde{\theta}_1,\ldots,\tilde{\theta}_d)$ is a constant vector, we observe from [2], page 600, and Lemma A.5 (see Appendix A) that

$$\frac{2^d\tilde{\theta}_1^3\cdots\tilde{\theta}_d^3}{\pi^d n^d}E\left[\tilde{X}_n'\left(\bigotimes_{t=1}^d R_{\tilde{\theta},n}\right)^{-1}\tilde{X}_n\right]$$

$$
\begin{aligned}
&= \frac{\phi^d \tilde{\theta}_1^3 \cdots \tilde{\theta}_d^3}{\theta_1^3 \cdots \theta_d^3 n^d} \operatorname{tr}\left[\left(\bigotimes_{t=1}^d R_{\tilde{\theta}_t,n}\right)^{-1}\left(\bigotimes_{t=1}^d R_{\theta_t,n}\right)\right] \\
&= \frac{\phi^d \tilde{\theta}_1^3 \cdots \tilde{\theta}_d^3}{\theta_1^3 \cdots \theta_d^3} \prod_{t=1}^d \left[\frac{1}{n} \operatorname{tr}(R_{\tilde{\theta}_t,n}^{-1} R_{\theta_t,n})\right] \\
&= \phi^d + O(n^{-1}),
\end{aligned}
$$

as $n \to \infty$. Next, we observe from [8], page 53, that

$$
\begin{aligned}
&\frac{2^{2d}\tilde{\theta}_1^6 \cdots \tilde{\theta}_d^6}{\pi^{2d} n^{2d}} E\left\{\left[\tilde{X}_n'\left(\bigotimes_{t=1}^d R_{\tilde{\theta}_t,n}\right)^{-1} \tilde{X}_n\right]^2\right\} \\
&\qquad = \frac{2\phi^{2d}\tilde{\theta}_1^6 \cdots \tilde{\theta}_d^6}{\theta_1^6 \cdots \theta_d^6 n^{2d}} \operatorname{tr}\left[\left(\bigotimes_{t=1}^d R_{\tilde{\theta}_t,n}^{-1} R_{\theta_t,n}\right)^2\right] \\
&\qquad\quad + \frac{\phi^{2d}\tilde{\theta}_1^6 \cdots \tilde{\theta}_d^6}{\theta_1^6 \cdots \theta_d^6 n^{2d}} \left[\operatorname{tr}\left(\bigotimes_{t=1}^d R_{\tilde{\theta}_t,n}^{-1} R_{\theta_t,n}\right)\right]^2.
\end{aligned}
$$

Hence, it follows from Lemma A.7 that

$$
\begin{aligned}
&\operatorname{Var}\left\{\frac{2^d \tilde{\theta}_1^3 \cdots \tilde{\theta}_d^3}{\pi^d n^d} \tilde{X}_n'\left(\bigotimes_{t=1}^d R_{\tilde{\theta}_t,n}\right)^{-1} \tilde{X}_n\right\} \\
&\quad = \frac{2\phi^{2d}\tilde{\theta}_1^6 \cdots \tilde{\theta}_d^6}{\theta_1^6 \cdots \theta_d^6 n^{2d}} \operatorname{tr}\left[\left(\bigotimes_{t=1}^d R_{\tilde{\theta}_t,n}^{-1} R_{\theta_t,n}\right)^2\right] \\
&\quad = \frac{2\phi^{2d}\tilde{\theta}_1^6 \cdots \tilde{\theta}_d^6}{\theta_1^6 \cdots \theta_d^6 n^d} \prod_{t=1}^d \left\{\frac{1}{n} \operatorname{tr}\left[\left(R_{\tilde{\theta}_t,n}^{-1} R_{\theta_t,n}\right)^2\right]\right\} \\
&\quad = O(n^{-d}),
\end{aligned}
\tag{27}
$$

as $n \to \infty$ uniformly over $\beta_{0,0} < \phi < \beta_{1,0}$ and $\beta_{0,t} < \theta_t, \tilde{\theta}_t < \beta_{1,t}, t = 1, \ldots, d$. The final statement of Theorem 9 follows from the definition of $\hat{\phi}_{\tilde{\theta}_1,\ldots,\tilde{\theta}_d}$, Chebyshev's inequality and the fact that the $d$th root function is continuous. □

Remark. Theorem 8 shows that the amount of Fisher information contained in the sample on the scale parameter $\phi$ is an order of magnitude greater than that on the correlation parameters $\theta_1, \ldots, \theta_d$. Thus, it should not be really surprising that incorrect specification of the values of $\theta_1, \ldots, \theta_d$ can still lead to consistent esimation of $\phi$. Crowder [4] (page 49) has a discussion of such a phenomenon in a more general setting (see also [12], page 175).

**7. Sieve maximum likelihood estimation.** Let $\tilde{X}_n \sim N_{n^d}(0, \Sigma_{\phi,\theta_1,\dots,\theta_d;n})$ be as in Section 1. The following definition is taken from [12], page 163.

DEFINITION. For a class of probability models $\{P_\xi : \xi \in \Xi\}$ for a random field on a given bounded domain, a function $h$ on $\Xi$ is said to be microergodic if and only if, for all $\xi$, $\xi' \in \Xi$, $h(\xi) \neq h(\xi')$ implies $P_\xi$ is orthogonal to $P_{\xi'}$.

We observe from [12], page 164, that, for $d = 1$, $\phi$ is microergodic, while $\theta_1$ is not. This is consistent with Theorem 8 which states that the Fisher information for $\theta_1$ is bounded as $n \to \infty$ if $d = 1$. We shall assume in this section that $d \geq 3$.

Let $0 < \nu < (d-2)/(d+1)$ and $0 < \beta_{0,t} < \beta_{1,t} < \infty$, $t = 0, \dots, d$, be known constants such that $\beta_{0,0} < \phi < \beta_{0,1}$ and $\beta_{0,t} < \theta_t < \beta_{1,t}$, $t = 1, \dots, d$. We define a sieve $\Omega_n$ on the parameter space of $(\phi, \theta_1, \dots, \theta_d)$, namely,

$$\Omega_n = \left\{ \left( \frac{i_0}{n^\nu}, \dots, \frac{i_d}{n^\nu} \right) : \beta_{0,t} \leq \frac{i_t}{n^\nu} \leq \beta_{1,t}, i_t \text{ integer}, 0 \leq t \leq d \right\}.$$

The sieve maximum likelihood estimator for $(\phi, \theta_1, \dots, \theta_d)$ is that element $(\hat{\phi}, \hat{\theta}_1, \dots, \hat{\theta}_d) \in \Omega_n$ such that

$$L_n(\hat{\phi}, \hat{\theta}_1, \dots, \hat{\theta}_d) = \sup\{L_n(\tilde{\phi}, \tilde{\theta}_1, \dots, \tilde{\theta}_d) : (\tilde{\phi}, \tilde{\theta}_1, \dots, \tilde{\theta}_d) \in \Omega_n\},$$

where the likelihood function $L_n$ is as in (5). For each $\varepsilon > 0$ sufficiently small, define

$$\Omega_{n,\varepsilon} = \Omega_n \backslash (\phi - \varepsilon, \phi + \varepsilon) \times (\theta_1 - \varepsilon, \theta_1 + \varepsilon) \times \cdots \times (\theta_d - \varepsilon, \theta_d + \varepsilon).$$

THEOREM 10. *Let $d \geq 3$, $\tilde{X}_n$ be as in Section 1 with covariance matrix $\Sigma_{\phi,\theta_1,\dots,\theta_d;n}$. Let $(\bar{\phi}, \bar{\theta}_1, \dots, \bar{\theta}_d)$ be a constant vector (depending only on $n$) such that $(\bar{\phi}, \bar{\theta}_1, \dots, \bar{\theta}_d) \in \Omega_n$, $|\phi - \bar{\phi}| < n^{-\nu}$ and $|\theta_t - \bar{\theta}_t| < n^{-\nu}$, $t = 1, \dots, d$. Then for any $\varepsilon > 0$,*

$$\lim_{n \to \infty} P\left\{ \frac{\sup\{L_n(\tilde{\phi}, \tilde{\theta}_1, \dots, \tilde{\theta}_d) : (\tilde{\phi}, \tilde{\theta}_1, \dots, \tilde{\theta}_d) \in \Omega_{n,\varepsilon}\}}{L_n(\bar{\phi}, \bar{\theta}_1, \dots, \bar{\theta}_d)} \leq \varepsilon \right\} = 1,$$

*where $L_n$ is the likelihood function given as in (5). Consequently, the sieve maximum likelihood estimate $(\hat{\phi}, \hat{\theta}_1, \dots, \hat{\theta}_d) \to (\phi, \theta_1, \dots, \theta_d)$ in probability as $n \to \infty$.*

PROOF. We observe from (23) that

$$\frac{2}{n^d}[\log L_n(\bar{\phi}, \bar{\theta}_1, \dots, \bar{\theta}_d) - \log L_n(\tilde{\phi}, \tilde{\theta}_1, \dots, \tilde{\theta}_d)]$$

$$= \log\left(\frac{\tilde{\phi}^d}{\bar{\phi}^d}\right) + \log\left(\frac{\bar{\theta}_1^3 \cdots \bar{\theta}_d^3}{\tilde{\theta}_1^3 \cdots \tilde{\theta}_d^3}\right)$$

$$+\frac{1}{n}\sum_{t=1}^{d}\log\biggl(\frac{|R_{\tilde{\theta}_t,n}|}{|R_{\bar{\theta}_t,n}|}\biggr)+\frac{2^d\tilde{\theta}_1^3\cdots\tilde{\theta}_d^3}{\pi^d\tilde{\phi}^d n^d}\tilde{X}_n'\Biggl(\bigotimes_{t=1}^{d}R_{\tilde{\theta}_t,n}\Biggr)^{-1}\tilde{X}_n$$
$$-\frac{2^d\bar{\theta}_1^3\cdots\bar{\theta}_d^3}{\pi^d\bar{\phi}^d n^d}\tilde{X}_n'\Biggl(\bigotimes_{t=1}^{d}R_{\bar{\theta}_t,n}\Biggr)^{-1}\tilde{X}_n.$$

Now,

$$P\Biggl\{\sup_{(\tilde{\phi},\tilde{\theta}_1,\ldots,\tilde{\theta}_d)\in\Omega_{n,\varepsilon}}\{\log[L_n(\tilde{\phi},\tilde{\theta}_1,\ldots,\tilde{\theta}_d)]-\log[L_n(\bar{\phi},\bar{\theta}_1,\ldots,\bar{\theta}_d)]\}\geq\log(\varepsilon)\Biggr\}$$
$$=P\Biggl\{\bigcup_{(\tilde{\phi},\tilde{\theta}_1,\ldots,\tilde{\theta}_d)\in\Omega_{n,\varepsilon}}\{\log[L_n(\tilde{\phi},\tilde{\theta}_1,\ldots,\tilde{\theta}_d)]$$
$$-\log[L_n(\bar{\phi},\bar{\theta}_1,\ldots,\bar{\theta}_d)]\}\geq\log(\varepsilon)\Biggr\}$$
$$\leq\sum_{(\tilde{\phi},\tilde{\theta}_1,\ldots,\tilde{\theta}_d)\in\Omega_{n,\varepsilon}}P\{\log[L_n(\bar{\phi},\bar{\theta}_1,\ldots,\bar{\theta}_d)]$$
$$-\log[L_n(\tilde{\phi},\tilde{\theta}_1,\ldots,\tilde{\theta}_d)]\leq-\log(\varepsilon)\}$$
$$=\sum_{(\tilde{\phi},\tilde{\theta}_1,\ldots,\tilde{\theta}_d)\in\Omega_{n,\varepsilon}}P\Biggl\{\frac{2^d\tilde{\theta}_1^3\cdots\tilde{\theta}_d^3}{\pi^d\tilde{\phi}^d n^d}\tilde{X}_n'\Biggl(\bigotimes_{t=1}^{d}R_{\tilde{\theta}_t,n}\Biggr)^{-1}\tilde{X}_n$$
$$-\frac{2^d\bar{\theta}_1^3\cdots\bar{\theta}_d^3}{\pi^d\bar{\phi}^d n^d}\tilde{X}_n'\Biggl(\bigotimes_{t=1}^{d}R_{\bar{\theta}_t,n}\Biggr)^{-1}\tilde{X}_n$$
$$-E\Biggl[\frac{2^d\tilde{\theta}_1^3\cdots\tilde{\theta}_d^3}{\pi^d\tilde{\phi}^d n^d}\tilde{X}_n'\Biggl(\bigotimes_{t=1}^{d}R_{\tilde{\theta}_t,n}\Biggr)^{-1}\tilde{X}_n$$
$$-\frac{2^d\bar{\theta}_1^3\cdots\bar{\theta}_d^3}{\pi^d\bar{\phi}^d n^d}\tilde{X}_n'\Biggl(\bigotimes_{t=1}^{d}R_{\bar{\theta}_t,n}\Biggr)^{-1}\tilde{X}_n\Biggr]$$
$$\leq-\log\biggl(\frac{\tilde{\phi}^d}{\bar{\phi}^d}\biggr)-\log\biggl(\frac{\bar{\theta}_1^3\cdots\bar{\theta}_d^3}{\tilde{\theta}_1^3\cdots\tilde{\theta}_d^3}\biggr)$$
$$-\frac{1}{n}\sum_{t=1}^{d}\log\biggl(\frac{|R_{\tilde{\theta}_t,n}|}{|R_{\bar{\theta}_t,n}|}\biggr)$$
$$-\frac{2^d\tilde{\theta}_1^3\cdots\tilde{\theta}_d^3}{\pi^d\tilde{\phi}^d n^d}E\tilde{X}_n'\Biggl(\bigotimes_{t=1}^{d}R_{\tilde{\theta}_t,n}\Biggr)^{-1}\tilde{X}_n$$

$$+ \frac{2^d \bar{\theta}_1^3 \cdots \bar{\theta}_d^3}{\pi^d \bar{\phi}^d n^d} E\tilde{X}_n' \left( \bigotimes_{t=1}^{d} R_{\bar{\theta}_t,n} \right)^{-1} \tilde{X}_n - \frac{2\log(\varepsilon)}{n^d} \biggr\}.$$

Hence, using Chebyshev's inequality, we obtain

$$\begin{aligned}
&P\biggl\{ \sup_{(\tilde{\phi},\tilde{\theta}_1,\ldots,\tilde{\theta}_d)\in\Omega_{n,\varepsilon}} \{\log[L_n(\tilde{\phi},\tilde{\theta}_1,\ldots,\tilde{\theta}_d)] - \log[L_n(\bar{\phi},\bar{\theta}_1,\ldots,\bar{\theta}_d)]\} \geq \log(\varepsilon) \biggr\} \\
&\quad \leq \sum_{(\tilde{\phi},\tilde{\theta}_1,\ldots,\tilde{\theta}_d)\in\Omega_{n,\varepsilon}} \mathrm{Var}\biggl\{ \frac{2^d \tilde{\theta}_1^3 \cdots \tilde{\theta}_d^3}{\pi^d \tilde{\phi}^d n^d} \tilde{X}_n' \left( \bigotimes_{t=1}^{d} R_{\tilde{\theta}_t,n} \right)^{-1} \tilde{X}_n \\
&\qquad\qquad - \frac{2^d \bar{\theta}_1^3 \cdots \bar{\theta}_d^3}{\pi^d \bar{\phi}^d n^d} \tilde{X}_n' \left( \bigotimes_{t=1}^{d} R_{\bar{\theta}_t,n} \right)^{-1} \tilde{X}_n \biggr\} \\
&\qquad \times \biggl\{ -\log\left( \frac{\tilde{\phi}^d}{\bar{\phi}^d} \right) - 3\log\left( \frac{\bar{\theta}_1 \cdots \bar{\theta}_d}{\tilde{\theta}_1 \cdots \tilde{\theta}_d} \right) \\
&\qquad\qquad - \frac{1}{n} \sum_{t=1}^{d} \log\left( \frac{|R_{\tilde{\theta}_t,n}|}{|R_{\bar{\theta}_t,n}|} \right) \qquad (28) \\
&\qquad\qquad - \frac{2^d \tilde{\theta}_1^3 \cdots \tilde{\theta}_d^3}{\pi^d \tilde{\phi}^d n^d} E\tilde{X}_n' \left( \bigotimes_{t=1}^{d} R_{\tilde{\theta}_t,n} \right)^{-1} \tilde{X}_n \\
&\qquad\qquad + \frac{2^d \bar{\theta}_1^3 \cdots \bar{\theta}_d^3}{\pi^d \bar{\phi}^d n^d} E\tilde{X}_n' \left( \bigotimes_{t=1}^{d} R_{\bar{\theta}_t,n} \right)^{-1} \tilde{X}_n - \frac{2\log(\varepsilon)}{n^d} \biggr\}^{-2}.
\end{aligned}$$

Next, we observe from Theorem 5 and Lemma A.5 (see Appendix A) that

$$\begin{aligned}
&\log\left( \frac{\tilde{\phi}^d}{\phi^d} \right) + 3\log\left( \frac{\theta_1 \cdots \theta_d}{\tilde{\theta}_1 \cdots \tilde{\theta}_d} \right) + \frac{1}{n} \sum_{t=1}^{d} \log\left( \frac{|R_{\tilde{\theta}_t,n}|}{|R_{\theta_t,n}|} \right) \\
&\quad + \frac{2^d \tilde{\theta}_1^3 \cdots \tilde{\theta}_d^3}{\pi^d \tilde{\phi}^d n^d} E\tilde{X}_n' \left( \bigotimes_{t=1}^{d} R_{\tilde{\theta}_t,n} \right)^{-1} \tilde{X}_n - 1 \\
&\quad = \log\left( \frac{\tilde{\phi}^d}{\phi^d} \right) + 3\log\left( \frac{\theta_1 \cdots \theta_d}{\tilde{\theta}_1 \cdots \tilde{\theta}_d} \right) \\
&\qquad + \frac{3n-4}{n} \sum_{t=1}^{d} \log\left( \frac{\tilde{\theta}_t}{\theta_t} \right) - \frac{2(n-2)}{n^2} \sum_{t=1}^{d} (\tilde{\theta}_t - \theta_t) \\
&\qquad + \frac{\phi^d \tilde{\theta}_1^3 \cdots \tilde{\theta}_d^3}{\tilde{\phi}^d \theta_1^3 \cdots \theta_d^3 n^d} \mathrm{tr}\left[ \bigotimes_{t=1}^{d} (R_{\tilde{\theta}_t,n}^{-1} R_{\theta_t,n}) \right] - 1 + O(n^{-2})
\end{aligned}$$

$$
\begin{aligned}
&= -\log\left(\frac{\phi^d}{\tilde{\phi}^d}\right) - 1 - \frac{4}{n}\sum_{t=1}^{d}\log\left(\frac{\tilde{\theta}_t}{\theta_t}\right) - \frac{2}{n}\sum_{t=1}^{d}(\tilde{\theta}_t - \theta_t) \\
&\quad + \left(\frac{\phi^d}{\tilde{\phi}^d}\right)\prod_{t=1}^{d}\left[\frac{\tilde{\theta}_t^3}{\theta_t^3 n}\operatorname{tr}(R_{\tilde{\theta}_t,n}^{-1}R_{\theta_t,n})\right] + O(n^{-2}) \\
&= -\log\left(\frac{\phi^d}{\tilde{\phi}^d}\right) - 1 - \frac{4}{n}\sum_{t=1}^{d}\log\left(\frac{\tilde{\theta}_t}{\theta_t}\right) - \frac{2}{n}\sum_{t=1}^{d}(\tilde{\theta}_t - \theta_t) \\
&\quad + \left(\frac{\phi^d}{\tilde{\phi}^d}\right)\prod_{t=1}^{d}\left\{\left(\frac{\tilde{\theta}_t}{\theta_t}\right)^3\left[\left(\frac{\theta_t}{\tilde{\theta}_t}\right)^3 + \frac{1}{n}\left[1 + \left(\frac{\theta_t}{\tilde{\theta}_t}\right)^2 - 2\left(\frac{\theta_t}{\tilde{\theta}_t}\right)^3\right.\right.\right. \\
&\qquad\qquad \left.\left.\left. - \frac{3\tilde{\theta}_t}{4}\left(\frac{\theta_t}{\tilde{\theta}_t}\right)^4 + \frac{\tilde{\theta}_t}{2}\left(\frac{\theta_t}{\tilde{\theta}_t}\right)^2 + \frac{\tilde{\theta}_t}{4}\right]\right]\right\} \\
&\quad + O(n^{-2}) \\
&= f\left(\frac{\phi^d}{\tilde{\phi}^d}\right) + \frac{\phi^d}{\tilde{\phi}^d n}\sum_{t=1}^{d} g(\theta_t, \tilde{\theta}_t) + O(n^{-2}),
\end{aligned}
$$

as $n \to \infty$ uniformly over $\beta_{0,0} \le \tilde{\phi} \le \beta_{0,1}$ and $\beta_{0,t} \le \tilde{\theta}_t \le \beta_{1,t}$, $t = 1, \ldots, d$, where

$$
\begin{aligned}
f(t) &= t - \log(t) - 1 \qquad \forall\, t > 0, \\
g(\theta_t, \tilde{\theta}_t) &= -4\log\left(\frac{\tilde{\theta}_t}{\theta_t}\right) - 2(\tilde{\theta}_t - \theta_t) + \left(\frac{\tilde{\theta}_t}{\theta_t}\right)^3 + \frac{\tilde{\theta}_t}{\theta_t} - 2 \\
&\quad - \frac{3\theta_t}{4} + \frac{\tilde{\theta}_t}{2}\left(\frac{\tilde{\theta}_t}{\theta_t}\right) + \frac{\tilde{\theta}_t}{4}\left(\frac{\tilde{\theta}_t}{\theta_t}\right)^3.
\end{aligned}
$$

We further observe that $f(t) \ge 0$ for all $t > 0$, $f(t) = 0$ if and only if $t = 1$,

$$
\begin{aligned}
\frac{\partial}{\partial \tilde{\theta}_t} g(\theta_t, \tilde{\theta}_t) &\ge 0 \qquad \forall\, \tilde{\theta}_t \ge \theta_t, \\
\frac{\partial}{\partial \tilde{\theta}_t} g(\theta_t, \tilde{\theta}_t) &\le 0 \qquad \forall\, \tilde{\theta}_t \le \theta_t, \\
g(\theta_t, \tilde{\theta}_t) &= 0 \qquad \text{if and only if } \theta_t = \tilde{\theta}_t.
\end{aligned}
$$

Hence, we conclude that

$$
\begin{aligned}
&\left\{ -\log\left(\frac{\tilde{\phi}^d}{\bar{\phi}^d}\right) - 3\log\left(\frac{\bar{\theta}_1 \cdots \bar{\theta}_d}{\tilde{\theta}_1 \cdots \tilde{\theta}_d}\right) - \frac{1}{n}\sum_{t=1}^{d}\log\left(\frac{|R_{\tilde{\theta}_t,n}|}{|R_{\bar{\theta}_t,n}|}\right)\right. \\
&\quad - \frac{2^d \tilde{\theta}_1^3 \cdots \tilde{\theta}_d^3}{\pi^d \tilde{\phi}^d n^d} E\tilde{X}_n'\left(\bigotimes_{t=1}^{d} R_{\tilde{\theta}_t,n}\right)^{-1}\tilde{X}_n
\end{aligned}
$$

$$
\begin{aligned}
& + \frac{2^d \bar\theta_1^3 \cdots \bar\theta_d^3}{\pi^d \bar\phi^d n^d} E \tilde X_n' \Bigg( \bigotimes_{t=1}^d R_{\bar\theta_t, n} \Bigg)^{-1} \tilde X_n - \frac{2\log(\varepsilon)}{n^d} \Bigg\}^{-2} \\
&= \Bigg\{ f\Bigg(\frac{\phi^d}{\tilde\phi^d}\Bigg) + \frac{\phi^d}{\tilde\phi^d n} \sum_{t=1}^d g(\theta_t, \tilde\theta_t) \\
&\qquad - f\Bigg(\frac{\phi^d}{\bar\phi^d}\Bigg) - \frac{\phi^d}{\bar\phi^d n} \sum_{t=1}^d g(\theta_t, \bar\theta_t) + O(n^{-2}) \Bigg\}^{-2} \\
&\le O(n^2),
\end{aligned} \tag{29}
$$

as $n \to \infty$ uniformly over $(\tilde\phi, \tilde\theta_1, \ldots, \tilde\theta_d) \in \Omega_{n,\varepsilon}$. Also, we observe from (27) that

$$
\begin{aligned}
& \operatorname{Var}\Bigg\{ \frac{2^d \tilde\theta_1^3 \cdots \tilde\theta_d^3}{\pi^d \tilde\phi^d n^d} \tilde X_n' \Bigg( \bigotimes_{t=1}^d R_{\tilde\theta_t, n} \Bigg)^{-1} \tilde X_n - \frac{2^d \bar\theta_1^3 \cdots \bar\theta_d^3}{\pi^d \bar\phi^d n^d} \tilde X_n' \Bigg( \bigotimes_{t=1}^d R_{\bar\theta_t, n} \Bigg)^{-1} \tilde X_n \Bigg\} \\
&\le 2 \operatorname{Var}\Bigg\{ \frac{2^d \tilde\theta_1^3 \cdots \tilde\theta_d^3}{\pi^d \tilde\phi^d n^d} \tilde X_n' \Bigg( \bigotimes_{t=1}^d R_{\tilde\theta_t, n} \Bigg)^{-1} \tilde X_n \Bigg\} \\
&\quad + 2 \operatorname{Var}\Bigg\{ \frac{2^d \bar\theta_1^3 \cdots \bar\theta_d^3}{\pi^d \bar\phi^d n^d} \tilde X_n' \Bigg( \bigotimes_{t=1}^d R_{\bar\theta_t, n} \Bigg)^{-1} \tilde X_n \Bigg\} \\
&= O(n^{-d}),
\end{aligned} \tag{30}
$$

as $n \to \infty$ uniformly over $\beta_{0,0} < \tilde\phi, \bar\phi < \beta_{1,0}$ and $\beta_{0,t} < \tilde\theta_t, \bar\theta_t < \beta_{1,t}, t = 1, \ldots, d$. Finally, since $0 < \nu < (d-2)/(d+1)$, it follows from (28), (29) and (30) that

$$
\begin{aligned}
& P\Bigg\{ \sup_{(\tilde\phi, \tilde\theta_1, \ldots, \tilde\theta_d) \in \Omega_{n,\varepsilon}} \{\log[L_n(\tilde\phi, \tilde\theta_1, \ldots, \tilde\theta_d)] - \log[L_n(\bar\phi, \bar\theta_1, \ldots, \bar\theta_d)]\} \ge \log(\varepsilon) \Bigg\} \\
&\quad = \sum_{(\tilde\phi, \tilde\theta_1, \ldots, \tilde\theta_d) \in \Omega_{n,\varepsilon}} O(n^{-d+2}) = O(n^{(d+1)\nu - d + 2}) \to 0
\end{aligned}
$$

as $n \to \infty$. This proves Theorem 10. □

## APPENDIX A

Lemma A.1. *With the notation of Theorem 1, the generating functions $G_\tau$ and $G_{\tau^*}$ of $\{\tau_k : k \ge 0\}$ and $\{\tau_k^* : k \ge 0\}$ are, respectively,*

$$
G_\tau(z) = \sum_{k=0}^\infty \frac{\tau_k z^k}{u^k} = \frac{a + bz}{(1-z)^2}, \qquad G_{\tau^*}(z) = \sum_{k=0}^\infty \frac{\tau_k^* z^k}{u^k} = \frac{\tilde a + \tilde b z}{(1-z)^2}.
$$

Proof. For all integers $k \ge 0$, we have

$$
\frac{\tau_k^*}{u^k} - \frac{\tau_{k-1}^*}{u^{k-1}} = (w - wu^2 - w^2u^2)\mathcal{I}\{k \ge 1\} + [1 - (1+w)^2u^2]\mathcal{I}\{k = 0\},
$$

where $\mathcal{I}\{\cdot\}$ denotes the indicator function. Hence

$$\sum_{k=0}^{\infty}\frac{\tau_k^* z^k}{u^k} - z\sum_{k=1}^{\infty}\frac{\tau_{k-1}^* z^{k-1}}{u^{k-1}} = (w - wu^2 - w^2u^2)\sum_{k=1}^{\infty} z^k + [1-(1+w)^2u^2],$$

which implies that

$$G_{\tau^*}(z) = \frac{1-(1+w)^2u^2}{1-z} + \frac{(w-wu^2-w^2u^2)z}{(1-z)^2}.$$

Similarly, for all integers $k \geq 0$,

$$\begin{aligned}\frac{\tau_k}{u^k} &- \frac{\tau_{k-1}}{u^{k-1}}\mathcal{I}\{k\geq 1\}\\ &= (w-2wu^2+wu^4)\mathcal{I}\{k\geq 1\} + [1-2(1+w)u^2+(1+2w)u^4]\mathcal{I}\{k=0\}.\end{aligned}$$

Hence

$$\sum_{k=0}^{\infty}\frac{\tau_k z^k}{u^k} - z\sum_{k=1}^{\infty}\frac{\tau_{k-1} z^{k-1}}{u^{k-1}} = (u^2-1)^2 + 2wu^2(u^2-1) + w(u^2-1)^2\sum_{k=1}^{\infty} z^k,$$

which implies that

$$G_\tau(z) = \frac{(u^2-1)^2+2wu^2(u^2-1)}{1-z} + \frac{w(u^2-1)^2 z}{(1-z)^2}. \qquad \square$$

LEMMA A.2. *For $1 \leq m \leq n-1$, let $\tilde{A}_{\theta,n;m}$ be the $m\times m$ matrix defined as in (11). Then with the notation of Theorem 1,*

$$\begin{aligned}|\tilde{A}_{\theta,n;m}| = (-1)^{m+1}\{&(a\alpha_1+u^2b)[(b+\tau_{-1}u)\alpha_1]^{m-1}\\ &- (a\alpha_2+u^2b)[(b+\tau_{-1}u)\alpha_2]^{m-1}\}/(\alpha_1-\alpha_2).\end{aligned}$$

PROOF. Using elementary row operations to reduce $\tilde{A}_{\theta,n;m}$ to a diagonal matrix, we obtain

$$\begin{aligned}|\tilde{A}_{\theta,n;m}| &= (-\tau_{-1})^{m-1}\tau_{m-1} + (-\tau_{-1})^{m-2}\sum_{k_1=1}^{m-1}\tau_{m-k_1-1}\tau_{k_1-1}\\ &\quad + (-\tau_{-1})^{m-3}\sum_{k_1=2}^{m-1}\sum_{k_2=1}^{k_1-1}\tau_{m-k_1-1}\tau_{k_1-k_2-1}\tau_{k_2-1}\\ &\quad + \cdots\\ &\quad - \tau_{-1}\sum_{k_1=m-2}^{m-1}\sum_{k_2=m-3}^{k_1-1}\cdots\sum_{k_{m-3}=2}^{k_{m-4}-1}\sum_{k_{m-2}=1}^{k_{m-3}-1}\tau_{m-k_1-1}\tau_{k_1-k_2-1}\cdots\\ &\qquad\qquad \times \tau_{k_{m-3}-k_{m-2}-1}\tau_{k_{m-2}-1}\\ &\quad + \tau_0^m.\end{aligned}$$

Letting $G_\tau(z)$ be as in Lemma A.1, we observe that $|\tilde{A}_{\theta,n;m}|$ equals the coefficient of $z^{m-1}$ in the series expansion of

$$\begin{aligned}
&(-\tau_{-1}u)^{m-1}G_\tau(z) + (-\tau_{-1}u)^{m-2}zG_\tau^2(z) \\
&\quad + (-\tau_{-1}u)^{m-3}z^2G_\tau^3(z) + \cdots - \tau_{-1}uz^{m-2}G_\tau^{m-1}(z) + z^{m-1}G_\tau^m(z) \\
&\quad = (-\tau_{-1}u)^{m-1}G_\tau(z)\sum_{k=0}^{m-1}(-\tau_{-1}u)^{-k}z^kG_\tau^k(z) \\
&\quad = (-\tau_{-1}u)^{m-1}G_\tau(z)\frac{(-\tau_{-1}u)^{-m}z^mG_\tau^m(z)-1}{(-\tau_{-1}u)^{-1}zG_\tau(z)-1}.
\end{aligned}$$

Following the proof of Theorem 1, we conclude that, for $1 \le m \le n-1$,

$$|\tilde{A}_{\theta,n;m}| = \frac{(-1)^{m+1}(b+\tau_{-1}u)^{m-1}[a(\alpha_1^m - \alpha_2^m) + u^2b(\alpha_1^{m-1} - \alpha_2^{m-1})]}{\alpha_1 - \alpha_2}.$$

This proves Lemma A.2. □

LEMMA A.3. *With the notation of Theorem 4, for $3 \le i \le j-2 \le n-2$,*

$$\begin{aligned}
&|R_{\theta,n;-i,-j;i+1}| \\
&\quad = \biggl\{\hat{\tau}_{-1}\tilde{\tau}_{-1} - \hat{\tau}_0\tilde{\tau}_{-2} + \frac{\hat{\tau}_{-2}u\tilde{\tau}_{-1}[w(1-4u^2+3u^4)+(1-u^2)^2]}{(b+\tau_{-1}u)\alpha_1 - \tau_{-1}u} \\
&\qquad - \frac{\hat{\tau}_{-2}\tilde{\tau}_{-2}u^2[w(2-6u^2+4u^4)+(1-u^2)^2]}{(b+\tau_{-1}u)\alpha_1 - \tau_{-1}u} \\
&\qquad\qquad + \frac{w\tau_{-1}u^2\hat{\tau}_{-2}(\tilde{\tau}_{-1} - u\tilde{\tau}_{-2})(1-u^2)^2}{[(b+\tau_{-1}u)\alpha_1 - \tau_{-1}u]^2}\biggr\} \\
&\qquad \times \frac{(-1)^{i-1}(\tilde{a}\alpha_1 + u^2\tilde{b})}{\alpha_1 - \alpha_2}[(b+\tau_{-1}u)\alpha_1]^{i-3} \\
&\qquad - \biggl\{\hat{\tau}_{-1}\tilde{\tau}_{-1} - \hat{\tau}_0\tilde{\tau}_{-2} + \frac{\hat{\tau}_{-2}u\tilde{\tau}_{-1}[w(1-4u^2+3u^4)+(1-u^2)^2]}{(b+\tau_{-1}u)\alpha_2 - \tau_{-1}u} \\
&\qquad - \frac{\hat{\tau}_{-2}\tilde{\tau}_{-2}u^2[w(2-6u^2+4u^4)+(1-u^2)^2]}{(b+\tau_{-1}u)\alpha_2 - \tau_{-1}u} \\
&\qquad\qquad + \frac{w\tau_{-1}u^2\hat{\tau}_{-2}(\tilde{\tau}_{-1} - u\tilde{\tau}_{-2})(1-u^2)^2}{[(b+\tau_{-1}u)\alpha_2 - \tau_{-1}u]^2}\biggr\} \\
&\qquad \times \frac{(-1)^{i-1}(\tilde{a}\alpha_2 + u^2\tilde{b})}{\alpha_1 - \alpha_2}[(b+\tau_{-1}u)\alpha_2]^{i-3}.
\end{aligned}$$

PROOF. First, for $l=1,2$ and $i-1\le m\le n-1$, define

$$S_{i,m,l}=\sum_{k=2}^{i-2}[(m-k)-2(m-k+1)u^2+(m-k+2)u^4] \\ \times wu^{m-k}\tau_{-1}^{-k}(b+\tau_{-1}u)^{k-2}\alpha_l^{k-2}.$$

Then

$$\begin{aligned}
&S_{i,m,l}-(b+\tau_{-1}u)\alpha_l u^{-1}\tau_{-1}^{-1}S_{i,m,l}\\
&\quad=\sum_{k=2}^{i-2}[(m-k)-2(m-k+1)u^2+(m-k+2)u^4]\\
&\qquad\times wu^{m-k}\tau_{-1}^{-k}(b+\tau_{-1}u)^{k-2}\alpha_l^{k-2}\\
&\qquad-\sum_{k=2}^{i-2}[(m-k)-2(m-k+1)u^2+(m-k+2)u^4]\\
&\qquad\times wu^{m-k-1}\tau_{-1}^{-k-1}(b+\tau_{-1}u)^{k-1}\alpha_l^{k-1}\\
&\quad=[(m-2)-2(m-1)u^2+mu^4]wu^{m-2}\tau_{-1}^{-2}\\
&\qquad-[(m-i+2)-2(m-i+3)u^2+(m-i+4)u^4]\\
&\qquad\times wu^{m-i+1}\tau_{-1}^{-i+1}(b+\tau_{-1}u)^{i-3}\alpha_l^{i-3}\\
&\qquad-\sum_{k=3}^{i-2}(1-2u^2+u^4)wu^{m-k}\tau_{-1}^{-k}(b+\tau_{-1}u)^{k-2}\alpha_l^{k-2}.
\end{aligned}$$

Hence, for $i\ge 3$,

$$\begin{aligned}
\frac{S_{i,m,l}}{u^m}&=\frac{w[(m-i+2)-2(m-i+3)u^2+(m-i+4)u^4]}{\tau_{-1}u[(b+\tau_{-1}u)\alpha_l-\tau_{-1}u]}\\
&\quad\times\left[\frac{(b+\tau_{-1}u)\alpha_l}{\tau_{-1}u}\right]^{i-3}\\
&\quad-\frac{w[(m-2)-2(m-1)u^2+mu^4]}{\tau_{-1}u[(b+\tau_{-1}u)\alpha_l-\tau_{-1}u]}\\
&\quad+\frac{w(1-u^2)^2}{[(b+\tau_{-1}u)\alpha_l-\tau_{-1}u]^2}\left[\frac{(b+\tau_{-1}u)\alpha_l}{\tau_{-1}u}\right]^{i-3}\\
&\quad-\frac{w(1-u^2)^2(b+\tau_{-1}u)\alpha_l}{\tau_{-1}u[(b+\tau_{-1}u)\alpha_l-\tau_{-1}u]^2}\\
&=\biggl\{\frac{w[(m-i+2)-2(m-i+3)u^2+(m-i+4)u^4]}{\tau_{-1}u[(b+\tau_{-1}u)\alpha_l-\tau_{-1}u]}
\end{aligned}\tag{31}$$

$$+ \frac{w(1-u^2)^2}{[(b+\tau_{-1}u)\alpha_l - \tau_{-1}u]^2}\Bigg\}$$

$$\times \left[\frac{(b+\tau_{-1}u)\alpha_l}{\tau_{-1}u}\right]^{i-3}$$

$$- \frac{w[(m-2)-2(m-1)u^2+mu^4]}{\tau_{-1}u[(b+\tau_{-1}u)\alpha_l - \tau_{-1}u]}$$

$$- \frac{w(1-u^2)^2(b+\tau_{-1}u)\alpha_l}{\tau_{-1}u[(b+\tau_{-1}u)\alpha_l - \tau_{-1}u]^2}.$$

Next we observe from Proposition 2 and Case 3 that, for $3 \le i \le j-2 \le n-2$,

$$\begin{aligned}
&|R_{\theta,n;-i,-j;i+1}| \\
&\quad = (-1)^i (A_{\theta,n;-i,-j;i+1})_{1,i+1} \prod_{l=0}^{i-1} (A_{\theta,n;-i,-j;i+1})_{i+1-l,i-l} \\
&\qquad + \sum_{k=1}^{i} (-1)^{i-k} (A_{\theta,n;-i,-j;i+1})_{k+1,i+1} |R_{\theta,n;-i,-j;k}| \\
&\qquad\qquad \times \prod_{l=0}^{i-k-1} (A_{\theta,n;-i,-j;i+1})_{i+1-l,i-l} \\
&\quad = (-1)^i (1+w)(1+iw) u^{i+1} \tilde{\tau}_{-2} \hat{\tau}_{-2} \tau_{-1}^{i-3} \qquad (32)\\
&\qquad + (-1)^{i-1} [1+(i-1)w] u^{i-1} \tilde{\tau}_{-2} \hat{\tau}_{-2} \tau_{-1}^{i-3} \\
&\qquad + \sum_{k=2}^{i-2} (-1)^{i-k} \tau_{i-k} |R_{\theta,n;-i,-j;k}| \tau_{-1}^{i-k-2} \hat{\tau}_{-2} \tilde{\tau}_{-2} \\
&\qquad - \hat{\tau}_0 |R_{\theta,n;-i,-j;i-1}| \tilde{\tau}_{-2} + \tilde{\tau}_{-1} |R_{\theta,n;-i,-j;i}| \\
&\quad = (-1)^{i-1} \tau_{i-1}^* \tau_{-1}^{i-3} \hat{\tau}_{-2} \tilde{\tau}_{-2} + \sum_{k=2}^{i-2} (-1)^{i-k} \tau_{i-k} |R_{\theta,n;k}| \tau_{-1}^{i-k-2} \hat{\tau}_{-2} \tilde{\tau}_{-2} \\
&\qquad - \hat{\tau}_0 |R_{\theta,n;i-1}| \tilde{\tau}_{-2} + \tilde{\tau}_{-1} |R_{\theta,n;-i,-i-1;i}|.
\end{aligned}$$

We further observe from Theorem 1 that, for $i-1 \le m \le n-1$,

$$\begin{aligned}
&\sum_{k=2}^{i-2} (-1)^{-k} \tau_{m-k} \tau_{-1}^{-k} |R_{\theta,n;k}| \\
&\quad = (1-u^2)^2 \sum_{k=2}^{i-2} (-1)^{-k} u^{m-k} \tau_{-1}^{-k} |R_{\theta,n;k}|
\end{aligned}$$

$$+\sum_{k=2}^{i-2}(-1)^{-k}[(m-k)-2(m-k+1)u^2+(m-k+2)u^4]$$
$$\times wu^{m-k}\tau_{-1}^{-k}|R_{\theta,n;k}|$$
$$=\frac{u^{m-1}(1-u^2)^2(\tilde{a}\alpha_1+u^2\tilde{b})}{\tau_{-1}(\alpha_1-\alpha_2)[\alpha_1(b+\tau_{-1}u)-\tau_{-1}u]}\left\{\left[\frac{\alpha_1(b+\tau_{-1}u)}{\tau_{-1}u}\right]^{i-3}-1\right\}$$
$$-\frac{u^{m-1}(1-u^2)^2(\tilde{a}\alpha_2+u^2\tilde{b})}{\tau_{-1}(\alpha_1-\alpha_2)[\alpha_2(b+\tau_{-1}u)-\tau_{-1}u]}\left\{\left[\frac{\alpha_2(b+\tau_{-1}u)}{\tau_{-1}u}\right]^{i-3}-1\right\}\tag{33}$$
$$+\frac{(\tilde{a}\alpha_1+u^2\tilde{b})S_{i,m,1}-(\tilde{a}\alpha_2+u^2\tilde{b})S_{i,m,2}}{\alpha_1-\alpha_2}.$$

It follows from (32) and (33) that

$$|R_{\theta,n;-i,-j;i+1}|=(-1)^{i-1}\tau_{i-1}^{*}\tau_{-1}^{i-3}\hat{\tau}_{-2}\tilde{\tau}_{-2}-\hat{\tau}_0|R_{\theta,n;i-1}|\tilde{\tau}_{-2}$$
$$+\tilde{\tau}_{-1}|R_{\theta,n;-i,-(i+1);i}|$$
$$+\frac{(-1)^i\tau_{-1}^{i-3}\hat{\tau}_{-2}\tilde{\tau}_{-2}u^{i-1}(1-u^2)^2(\tilde{a}\alpha_1+u^2\tilde{b})}{(\alpha_1-\alpha_2)[\alpha_1(b+\tau_{-1}u)-\tau_{-1}u]}$$
$$\times\left\{\left[\frac{\alpha_1(b+\tau_{-1}u)}{\tau_{-1}u}\right]^{i-3}-1\right\}$$
$$-\frac{(-1)^i\tau_{-1}^{i-3}\hat{\tau}_{-2}\tilde{\tau}_{-2}u^{i-1}(1-u^2)^2(\tilde{a}\alpha_2+u^2\tilde{b})}{(\alpha_1-\alpha_2)[\alpha_2(b+\tau_{-1}u)-\tau_{-1}u]}$$
$$\times\left\{\left[\frac{\alpha_2(b+\tau_{-1}u)}{\tau_{-1}u}\right]^{i-3}-1\right\}$$
$$+\frac{(-1)^i\tau_{-1}^{i-2}\hat{\tau}_{-2}\tilde{\tau}_{-2}[(\tilde{a}\alpha_1+u^2\tilde{b})S_{i,i,1}-(\tilde{a}\alpha_2+u^2\tilde{b})S_{i,i,2}]}{\alpha_1-\alpha_2},$$

and Lemma A.3 results using (14) after some algebra. □

LEMMA A.4. *Let $\xi_k$ be defined as in (17). Then the generating function $G_\xi$ of $\{\xi_k : k\ge i+1\}$ is given by*

$$G_\xi(z)=\sum_{k=i+1}^{\infty}\frac{\xi_k z^k}{u^k}=\frac{z^{i+1}(a_\xi+b_\xi z)}{(1-z)^2},$$

*where*

$$a_\xi=\left\{\frac{\hat{\tau}_{-2}\tilde{\tau}_0[w(2-6u^2+4u^4)+(1-u^2)^2]}{\tau_{-1}u^2[(b+\tau_{-1}u)\alpha_1-\tau_{-1}u]}\right.$$
$$+\frac{w\hat{\tau}_{-2}\tilde{\tau}_0(1-u^2)^2}{u[(b+\tau_{-1}u)\alpha_1-\tau_{-1}u]^2}+\frac{\hat{\tau}_0\tilde{\tau}_0-\hat{\tau}_1\tilde{\tau}_{-1}}{\tau_{-1}u^4}$$

$$
\begin{aligned}
&- \frac{\hat\tau_{-2}\tilde\tau_{-1}[w(3-8u^2+5u^4)+(1-u^2)^2]}{\tau_{-1}u[(b+\tau_{-1}u)\alpha_1-\tau_{-1}u]} - \frac{w\hat\tau_{-2}\tilde\tau_{-1}(1-u^2)^2}{[(b+\tau_{-1}u)\alpha_1-\tau_{-1}u]^2}\Big\} \\
&\times \frac{(\tilde a\alpha_1+u^2\tilde b)}{\tau_{-1}^3(\alpha_1-\alpha_2)}\Big[\frac{(b+\tau_{-1}u)\alpha_1}{\tau_{-1}u}\Big]^{i-3} \\
&- \Big\{\frac{\hat\tau_{-2}\tilde\tau_0[w(2-6u^2+4u^4)+(1-u^2)^2]}{\tau_{-1}u^2[(b+\tau_{-1}u)\alpha_2-\tau_{-1}u]} \\
&\quad + \frac{w\hat\tau_{-2}\tilde\tau_0(1-u^2)^2}{u[(b+\tau_{-1}u)\alpha_2-\tau_{-1}u]^2} + \frac{\hat\tau_0\tilde\tau_0-\hat\tau_1\tilde\tau_{-1}}{\tau_{-1}u^4} \\
&\quad - \frac{\hat\tau_{-2}\tilde\tau_{-1}[w(3-8u^2+5u^4)+(1-u^2)^2]}{\tau_{-1}u[(b+\tau_{-1}u)\alpha_2-\tau_{-1}u]} - \frac{w\hat\tau_{-2}\tilde\tau_{-1}(1-u^2)^2}{[(b+\tau_{-1}u)\alpha_2-\tau_{-1}u]^2}\Big\} \\
&\times \frac{(\tilde a\alpha_2+u^2\tilde b)}{\tau_{-1}^3(\alpha_1-\alpha_2)}\Big[\frac{(b+\tau_{-1}u)\alpha_2}{\tau_{-1}u}\Big]^{i-3}
\end{aligned}
$$

*and*

$$
\begin{aligned}
a_\xi + b_\xi = {}& \Big\{\frac{w\hat\tau_{-2}u^2(2-3u^2+u^6)[w(2-6u^2+4u^4)+(1-u^2)^2]}{2\tau_{-1}u^4[(b+\tau_{-1}u)\alpha_1-\tau_{-1}u]} \\
&+ \frac{w^2\hat\tau_{-2}(2-3u^2+u^6)(1-u^2)^2}{2u[(b+\tau_{-1}u)\alpha_1-\tau_{-1}u]^2} + \frac{w\hat\tau_0(2-3u^2+u^6)}{2\tau_{-1}u^4} \\
&- \frac{w(1-u^2)^2\hat\tau_{-2}\tilde\tau_{-1}}{\tau_{-1}u[(b+\tau_{-1}u)\alpha_1-\tau_{-1}u]} - \frac{w(1-3u^4+2u^6)\tilde\tau_{-1}}{\tau_{-1}u^3}\Big\} \\
&\times \frac{(\tilde a\alpha_1+u^2\tilde b)}{\tau_{-1}^3(\alpha_1-\alpha_2)}\Big[\frac{(b+\tau_{-1}u)\alpha_1}{\tau_{-1}u}\Big]^{i-3} \\
&- \Big\{\frac{w\hat\tau_{-2}u^2(2-3u^2+u^6)[w(2-6u^2+4u^4)+(1-u^2)^2]}{2\tau_{-1}u^4[(b+\tau_{-1}u)\alpha_2-\tau_{-1}u]} \\
&\quad + \frac{w^2\hat\tau_{-2}(2-3u^2+u^6)(1-u^2)^2}{2u[(b+\tau_{-1}u)\alpha_2-\tau_{-1}u]^2} + \frac{w\hat\tau_0(2-3u^2+u^6)}{2\tau_{-1}u^4} \\
&\quad - \frac{w(1-u^2)^2\hat\tau_{-2}\tilde\tau_{-1}}{\tau_{-1}u[(b+\tau_{-1}u)\alpha_2-\tau_{-1}u]} - \frac{w(1-3u^4+2u^6)\tilde\tau_{-1}}{\tau_{-1}u^3}\Big\} \\
&\times \frac{(\tilde a\alpha_2+u^2\tilde b)}{\tau_{-1}^3(\alpha_1-\alpha_2)}\Big[\frac{(b+\tau_{-1}u)\alpha_2}{\tau_{-1}u}\Big]^{i-3}.
\end{aligned}
$$

*Next, let $\tilde{\xi}_k$ be defined as in (20). Then the generating function $G_{\tilde{\xi}}$ of $\{\tilde{\xi}_k : k \geq i+1\}$ is given by*

$$G_{\tilde{\xi}}(z) = \sum_{k=i+1}^{\infty} \frac{\tilde{\xi}_k z^k}{u^k} = \frac{z^{i+1}(a_{\tilde{\xi}} + b_{\tilde{\xi}} z)}{(1-z)^2},$$

*where*

$$\begin{aligned} a_{\tilde{\xi}} = \biggl\{ & \frac{\hat{\tau}_{-1}\tilde{\tau}_0}{\tau_{-1}u^4} + \frac{\hat{\tau}_{-2}\tilde{\tau}_0[w(1-4u^2+3u^4)+(1-u^2)^2]}{\tau_{-1}u^3[(b+\tau_{-1}u)\alpha_1 - \tau_{-1}u]} \\ & + \frac{w\tau_{-1}\hat{\tau}_{-2}\tilde{\tau}_0(1-u^2)^2}{\tau_{-1}u^2[(b+\tau_{-1}u)\alpha_1 - \tau_{-1}u]^2} \\ & - \frac{\hat{\tau}_1\tilde{\tau}_{-2}}{\tau_{-1}u^4} - \frac{\hat{\tau}_{-2}\tilde{\tau}_{-2}[w(3-8u^2+5u^4)+(1-u^2)^2]}{\tau_{-1}u[(b+\tau_{-1}u)\alpha_1 - \tau_{-1}u]} \\ & - \frac{w\hat{\tau}_{-2}\tilde{\tau}_{-2}(1-u^2)^2}{[(b+\tau_{-1}u)\alpha_1 - \tau_{-1}u]^2} \biggr\} \\ & \times \frac{\tilde{a}\alpha_1 + u^2\tilde{b}}{\tau_{-1}^3(\alpha_1-\alpha_2)} \left[\frac{(b+\tau_{-1}u)\alpha_1}{\tau_{-1}u}\right]^{i-3} \\ & - \biggl\{ \frac{\hat{\tau}_{-1}\tilde{\tau}_0}{\tau_{-1}u^4} + \frac{\hat{\tau}_{-2}\tilde{\tau}_0[w(1-4u^2+3u^4)+(1-u^2)^2]}{\tau_{-1}u^3[(b+\tau_{-1}u)\alpha_2 - \tau_{-1}u]} \\ & + \frac{w\tau_{-1}\hat{\tau}_{-2}\tilde{\tau}_0(1-u^2)^2}{\tau_{-1}u^2[(b+\tau_{-1}u)\alpha_2 - \tau_{-1}u]^2} \\ & - \frac{\hat{\tau}_1\tilde{\tau}_{-2}}{\tau_{-1}u^4} - \frac{\hat{\tau}_{-2}\tilde{\tau}_{-2}[w(3-8u^2+5u^4)+(1-u^2)^2]}{\tau_{-1}u[(b+\tau_{-1}u)\alpha_2 - \tau_{-1}u]} \\ & - \frac{w\hat{\tau}_{-2}\tilde{\tau}_{-2}(1-u^2)^2}{[(b+\tau_{-1}u)\alpha_2 - \tau_{-1}u]^2} \biggr\} \\ & \times \frac{\tilde{a}\alpha_2 + u^2\tilde{b}}{\tau_{-1}^3(\alpha_1-\alpha_2)} \left[\frac{(b+\tau_{-1}u)\alpha_2}{\tau_{-1}u}\right]^{i-3} \end{aligned}$$

*and*

$$\begin{aligned} a_{\tilde{\xi}} + b_{\tilde{\xi}} = \biggl\{ & \frac{w\hat{\tau}_{-1}(2-3u^2+u^6)}{2\tau_{-1}u^4} \\ & + \frac{w\hat{\tau}_{-2}(2-3u^2+u^6)[w(1-4u^2+3u^4)+(1-u^2)^2]}{2\tau_{-1}u^3[(b+\tau_{-1}u)\alpha_1 - \tau_{-1}u]} \\ & + \frac{w^2\hat{\tau}_{-2}(2-3u^2+u^6)(1-u^2)^2}{2u^2[(b+\tau_{-1}u)\alpha_1 - \tau_{-1}u]^2} - \frac{w\hat{\tau}_{-2}\tilde{\tau}_{-2}(1-u^2)^2}{\tau_{-1}u[(b+\tau_{-1}u)\alpha_1 - \tau_{-1}u]} \end{aligned}$$

$$
\begin{aligned}
&\qquad\qquad\qquad\qquad\qquad\qquad\qquad\qquad - \frac{w(1-3u^4+2u^6)\tilde{\tau}_{-2}}{\tau_{-1}u^3}\Bigg\} \\
&\times \frac{\tilde{a}\alpha_1 + u^2\tilde{b}}{\tau_{-1}^3(\alpha_1-\alpha_2)}\left[\frac{(b+\tau_{-1}u)\alpha_1}{\tau_{-1}u}\right]^{i-3} \\
&- \Bigg\{\frac{w\hat{\tau}_{-1}(2-3u^2+u^6)}{2\tau_{-1}u^4} \\
&\quad + \frac{w\hat{\tau}_{-2}(2-3u^2+u^6)[w(1-4u^2+3u^4)+(1-u^2)^2]}{2\tau_{-1}u^3[(b+\tau_{-1}u)\alpha_2-\tau_{-1}u]} \\
&\quad + \frac{w^2\hat{\tau}_{-2}(2-3u^2+u^6)(1-u^2)^2}{2u^2[(b+\tau_{-1}u)\alpha_2-\tau_{-1}u]^2} - \frac{w\hat{\tau}_{-2}\tilde{\tau}_{-2}(1-u^2)^2}{\tau_{-1}u[(b+\tau_{-1}u)\alpha_2-\tau_{-1}u]} \\
&\qquad\qquad\qquad\qquad\qquad\qquad\qquad\qquad - \frac{w(1-3u^4+2u^6)\tilde{\tau}_{-2}}{\tau_{-1}u^3}\Bigg\} \\
&\times \frac{\tilde{a}\alpha_2 + u^2\tilde{b}}{\tau_{-1}^3(\alpha_1-\alpha_2)}\left[\frac{(b+\tau_{-1}u)\alpha_2}{\tau_{-1}u}\right]^{i-3}.
\end{aligned}
$$

PROOF. First we observe after some algebra that $a_\xi = \xi_{i+1}/u^{i+1}$ and $a_\xi + b_\xi = \xi_k u^{-k} - \xi_{k-1}u^{-k+1}$, $k \geq i+2$. Hence, for integers $k \geq i+1$, we have

$$\frac{\xi_k}{u^k} - \frac{\xi_{k-1}}{u^{k-1}}\mathcal{I}\{k \geq i+2\} = a_\xi \mathcal{I}\{k = i+1\} + (a_\xi + b_\xi)\mathcal{I}\{k \geq i+2\}.$$

Hence

$$\sum_{k=i+1}^{\infty} \frac{\xi_k z^k}{u^k} - z\sum_{k=i+2}^{\infty} \frac{\xi_{k-1}z^{k-1}}{u^{k-1}} = a_\xi z^{i+1} + (a_\xi + b_\xi)\sum_{k=i+2}^{\infty} z^k$$

and

$$G_\xi(z) = \frac{z^{i+1}(a_\xi + b_\xi z)}{(1-z)^2}.$$

Similarly, we observe after some algebra that $a_{\tilde{\xi}} = \tilde{\xi}_{i+1}/u^{i+1}$ and $a_{\tilde{\xi}} + b_{\tilde{\xi}} = \tilde{\xi}_k u^{-k} - \tilde{\xi}_{k-1}u^{-k+1}$, $k \geq i+2$. Hence, for integers $k \geq i+1$, we have

$$\frac{\tilde{\xi}_k}{u^k} - \frac{\tilde{\xi}_{k-1}}{u^{k-1}}\mathcal{I}\{k \geq i+2\} = a_{\tilde{\xi}} \mathcal{I}\{k = i+1\} + (a_{\tilde{\xi}} + b_{\tilde{\xi}})\mathcal{I}\{k \geq i+2\}.$$

Hence

$$\sum_{k=i+1}^{\infty} \frac{\tilde{\xi}_k z^k}{u^k} - z\sum_{k=i+2}^{\infty} \frac{\tilde{\xi}_{k-1}z^{k-1}}{u^{k-1}} = a_{\tilde{\xi}} z^{i+1} + (a_{\tilde{\xi}} + b_{\tilde{\xi}})\sum_{k=i+2}^{\infty} z^k$$

and

$$G_{\tilde{\xi}}(z) = \frac{z^{i+1}(a_{\tilde{\xi}} + b_{\tilde{\xi}} z)}{(1-z)^2}.$$

This proves Lemma A.4. □

LEMMA A.5. *Let $0 < \beta_0 \le \theta, \tilde{\theta} \le \beta_1 < \infty$. Then*

$$\frac{1}{n} \operatorname{tr}(R_{\theta,n}^{-1} R_{\tilde{\theta},n}) = \left(\frac{\tilde{\theta}}{\theta}\right)^3 - \frac{w}{4}\left[3\left(\frac{\tilde{\theta}}{\theta}\right)^4 - 2\left(\frac{\tilde{\theta}}{\theta}\right)^2 - 1\right] + \frac{1}{n}\left[1 + \left(\frac{\tilde{\theta}}{\theta}\right)^2 - 2\left(\frac{\tilde{\theta}}{\theta}\right)^3\right] + O(n^{-2}),$$

*as $n \to \infty$ uniformly over $\beta_0 \le \theta, \tilde{\theta} \le \beta_1$.*

PROOF. Let $\tilde{w} = \tilde{\theta}/n$ and and $m_n = \lfloor (n+3)/2 \rfloor$, where $\lfloor \cdot \rfloor$ denotes the greatest integer function. We further define

$$r_{\alpha_1} = \frac{e^{-\tilde{w}} \tau_{-1}}{(b + \tau_{-1} u)\alpha_1} = \frac{e^{-\tilde{w}} u}{\alpha_1} = -\frac{1}{2+\sqrt{3}}\left[1 + O\left(\frac{1}{n}\right)\right]$$

and

$$r_{\alpha_2} = \frac{e^{-\tilde{w}} \tau_{-1}}{(b + \tau_{-1} u)\alpha_2} = \frac{e^{-\tilde{w}} u}{\alpha_2} = -\frac{1}{2-\sqrt{3}}\left[1 + O\left(\frac{1}{n}\right)\right],$$

as $n \to \infty$ uniformly over $\beta_0 \le \theta, \tilde{\theta} \le \beta_1$. Using *Mathematica*, we observe from Theorem 6 that, as $n \to \infty$,

$$\begin{aligned}(R_{\tilde{\theta},n} R_{\theta,n}^{-1})_{n,n} &= (R_{\tilde{\theta},n} R_{\theta,n}^{-1})_{1,1} \\ &= \sum_{i=1}^{n} [1 + (i-1)\tilde{w}] e^{-(i-1)\tilde{w}} (R_{\theta,n}^{-1})_{i,1} \\ &= \frac{3(70226 + 40545\sqrt{3})}{(664626 + 383722\sqrt{3})w}\left[1 - \left(\frac{\tilde{\theta}}{\theta}\right)^2\right] \\ &\quad + \frac{1713 + 989\sqrt{3}}{4(26 + 15\sqrt{3})^2} + \frac{989 + 571\sqrt{3}}{4(26 + 15\sqrt{3})^2}\left(\frac{\tilde{\theta}}{\theta}\right)^2 \\ &\quad + \frac{2702 + 1560\sqrt{3}}{4(26 + 15\sqrt{3})^2}\left(\frac{\tilde{\theta}}{\theta}\right)^3 + O(n^{-1}),\end{aligned} \tag{34}$$

as $n \to \infty$ uniformly over $\beta_0 \le \theta, \tilde{\theta} \le \beta_1$. Also, we observe from Theorem 7 that

$$(R_{\tilde{\theta},n} R_{\theta,n}^{-1})_{i,i}$$

$$
\begin{aligned}
&= \sum_{j=1}^{n} (1 + |i-j|\tilde{w}) e^{-|i-j|\tilde{w}} (R_{\theta,n}^{-1})_{j,i} \\
&= [1 + (i-1)\tilde{w}] e^{-(i-1)\tilde{w}} (R_{\theta,n}^{-1})_{n-i+1,n} \\
&\quad + [1 + (i-2)\tilde{w}] e^{-(i-2)\tilde{w}} (R_{\theta,n}^{-1})_{n-i+1,n-1} \\
&\quad + (1+\tilde{w}) e^{-\tilde{w}} (R_{\theta,n}^{-1})_{i-1,i} + (R_{\theta,n}^{-1})_{i,i} \\
&\quad + (1+\tilde{w}) e^{-\tilde{w}} (R_{\theta,n}^{-1})_{i,i+1} + [1 + (n-i)\tilde{w}] e^{-(n-i)\tilde{w}} (R_{\theta,n}^{-1})_{i,n} \\
&\quad + \biggl\{ \frac{C_{6,1} u^2}{\tau_{-1}^4 \alpha_1^2 (\alpha_1 - \alpha_2)} \biggl[ (a\alpha_1 + u^2 b) - (a\alpha_2 + u^2 b) \biggl(\frac{\alpha_2}{\alpha_1}\biggr)^{n-i-1} \biggr] \\
&\qquad \times \{ r_{\alpha_1}^4 (1 + 2\tilde{w}) - r_{\alpha_1}^5 (1+\tilde{w}) \qquad (35) \\
&\qquad + r_{\alpha_1}^{i+1} [1 + (i-3)\tilde{w}] - r_{\alpha_1}^{i} [1 + (i-2)\tilde{w}] \} / \{ r_{\alpha_1}^2 (1 - r_{\alpha_1})^2 \} \\
&\quad - \frac{C_{6,2} u^2 r_{\alpha_1}^i (\tilde{a}\alpha_2 + u^2 \tilde{b})}{\tau_{-1}^4 \alpha_1^2 (\alpha_1 - \alpha_2)(\tilde{a}\alpha_1 + u^2 \tilde{b})} \\
&\qquad \times \biggl[ (a\alpha_1 + u^2 b) \biggl(\frac{\alpha_2}{\alpha_1}\biggr)^{-3} - (a\alpha_2 + u^2 b) \biggl(\frac{\alpha_2}{\alpha_1}\biggr)^{n-i-4} \biggr] \\
&\qquad \times \{ r_{\alpha_2}^{-1} [1 + (i-3)\tilde{w}] - r_{\alpha_2}^{-2} [1 + (i-2)\tilde{w}] \\
&\qquad\qquad - r_{\alpha_2}^{-i+3} (1+\tilde{w}) + r_{\alpha_2}^{-i+2} (1 + 2\tilde{w}) \} / \{ (1 - r_{\alpha_2})^2 \} \\
&\quad + \frac{C_{6,1} u^2 r_{\alpha_1}^{-i} (a\alpha_1 + u^2 b)}{\tau_{-1}^4 \alpha_1^2 (\alpha_1 - \alpha_2)} \\
&\qquad \times \{ r_{\alpha_1}^{i+2} (1 + 2\tilde{w}) \\
&\qquad\quad - r_{\alpha_1}^{i+3} (1+\tilde{w}) - r_{\alpha_1}^{n} [1 + (n-i)\tilde{w}] \\
&\qquad\qquad + r_{\alpha_1}^{n+1} [1 + (n-i-1)\tilde{w}] \} / \{ (1 - r_{\alpha_1})^2 \} \\
&\quad - \frac{C_{6,1} u^2 r_{\alpha_1}^{-i} (a\alpha_2 + u^2 b)}{\tau_{-1}^4 \alpha_1^2 (\alpha_1 - \alpha_2)} \biggl(\frac{\alpha_2}{\alpha_1}\biggr)^{n-1} \\
&\qquad \times \{ r_{\alpha_2}^{i+2} (1 + 2\tilde{w}) - r_{\alpha_2}^{i+3} (1+\tilde{w}) \\
&\qquad\quad - r_{\alpha_2}^{n} [1 + (n-i)\tilde{w}] \\
&\qquad\quad + r_{\alpha_2}^{n+1} [1 + (n-i-1)\tilde{w}] \} / \{ (1 - r_{\alpha_2})^2 \} \\
&\quad - \frac{C_{6,2} u^2 r_{\alpha_1}^{-i} (a\alpha_1 + u^2 b)(\tilde{a}\alpha_2 + u^2 \tilde{b})}{\tau_{-1}^4 \alpha_1^2 (\alpha_1 - \alpha_2)(\tilde{a}\alpha_1 + u^2 \tilde{b})} \biggl(\frac{\alpha_2}{\alpha_1}\biggr)^{i-3} \\
&\qquad \times \{ r_{\alpha_1}^{i+2} (1 + 2\tilde{w}) - r_{\alpha_1}^{i+3} (1+\tilde{w})
\end{aligned}
$$

$$
\begin{aligned}
&\qquad - r_{\alpha_1}^{n}[1+(n-i)\tilde{w}] \\
&\qquad\quad + r_{\alpha_1}^{n+1}[1+(n-i-1)\tilde{w}]\}/\{(1-r_{\alpha_1})^2\} \\
&+ \frac{C_{6,2}u^2 r_{\alpha_1}^{-i}(a\alpha_2+u^2b)(\tilde{a}\alpha_2+u^2\tilde{b})}{\tau_{-1}^4\alpha_1^2(\alpha_1-\alpha_2)(\tilde{a}\alpha_1+u^2\tilde{b})}\left(\frac{\alpha_2}{\alpha_1}\right)^{n+i-4} \\
&\quad \times \{r_{\alpha_2}^{i+2}(1+2\tilde{w}) - r_{\alpha_2}^{i+3}(1+\tilde{w}) \\
&\qquad - r_{\alpha_2}^{n}[1+(n-i)\tilde{w}] + r_{\alpha_2}^{n+1}[1+(n-i-1)\tilde{w}]\}/\{(1-r_{\alpha_2})^2\}\Big\} \\
&\times \left[1+O\left(\left(\frac{\alpha_2}{\alpha_1}\right)^{n}\right)\right],
\end{aligned}
$$

as $n \to \infty$ uniformly over $5 \le i \le n-2$ and $\beta_0 \le \theta, \tilde{\theta} \le \beta_1$. Consequently, it follows from (35) that

$$
\begin{aligned}
&\sum_{i=m_n}^{n-2} (R_{\tilde{\theta},n}R_{\theta,n}^{-1})_{i,i} \\
&\quad = (n-m_n-1)\left\{\left(\frac{\tilde{\theta}}{\theta}\right)^3 - \frac{w}{4}\left[3\left(\frac{\tilde{\theta}}{\theta}\right)^4 - 2\left(\frac{\tilde{\theta}}{\theta}\right)^2 - 1\right] + O(n^{-2})\right\} \\
&\qquad + \frac{3(36810643322+21252634831\sqrt{3})}{(1300170624726+750653860178\sqrt{3})w}\left[1-\left(\frac{\tilde{\theta}}{\theta}\right)^2\right] \\
&\qquad - \frac{6322680+3650401\sqrt{3}}{128934018+74440090\sqrt{3}}\left[1-\left(\frac{\tilde{\theta}}{\theta}\right)^2\right] \\
&\qquad - \frac{39176289+22618441\sqrt{3}}{128934018+74440090\sqrt{3}}\left(1-\frac{\tilde{\theta}}{\theta}\right)\left(\frac{\tilde{\theta}}{\theta}\right)^2 + O(n^{-1}),
\end{aligned} \tag{36}
$$

as $n \to \infty$ uniformly over $\beta_0 \le \theta, \tilde{\theta} \le \beta_1$. Furthermore, if $n$ is an odd integer, we have

$$
\begin{aligned}
&(R_{\tilde{\theta},n}R_{\theta,n}^{-1})_{(n+1)/2,(n+1)/2} \\
&\quad = \left(\frac{\tilde{\theta}}{\theta}\right)^3 - \left[3\left(\frac{\tilde{\theta}}{\theta}\right)^4 - 2\left(\frac{\tilde{\theta}}{\theta}\right)^2 - 1\right]\frac{w}{4} + O(n^{-2}),
\end{aligned} \tag{37}
$$

as $n \to \infty$ uniformly over $\beta_0 \le \theta, \tilde{\theta} \le \beta_1$. Next we observe from Theorems 6 and 7 that

$$
\begin{aligned}
&(R_{\tilde{\theta},n}R_{\theta,n}^{-1})_{n-1,n-1} \\
&\quad = \sum_{j=1}^{n}(1+|n-1-j|\tilde{w})e^{-|n-1-j|\tilde{w}}(R_{\theta,n}^{-1})_{j,n-1}
\end{aligned}
$$

$$
\begin{aligned}
&= \frac{3(3351044259 + 1934726305\sqrt{3}\,)}{(25012534866 + 14440993738\sqrt{3}\,)w}\left[\left(\frac{\tilde{\theta}}{\theta}\right)^2 - 1\right] \\
&\quad + \frac{8016837 + 4628523\sqrt{3}}{34547766 + 19946162\sqrt{3}} \\
&\quad + \frac{3(6584767 + 3801717\sqrt{3})}{34547766 + 19946162\sqrt{3}}\left(\frac{\tilde{\theta}}{\theta}\right)^2 \\
&\quad + \frac{4(1694157 + 978122\sqrt{3})}{34547766 + 19946162\sqrt{3}}\left(\frac{\tilde{\theta}}{\theta}\right)^3 + O(n^{-1}),
\end{aligned} \tag{38}
$$

as $n \to \infty$ uniformly over $\beta_0 \le \theta, \tilde{\theta} \le \beta_1$. Finally, we conclude from (34), (36), (37) and (38) that

$$
\begin{aligned}
\frac{1}{n}\sum_{i=1}^{n}(R_{\tilde{\theta},n}R_{\theta,n}^{-1})_{i,i} &= \frac{2}{n}\sum_{i=m_n}^{n}(R_{\tilde{\theta},n}R_{\theta,n}^{-1})_{i,i} + \frac{1}{n}\left(\frac{\tilde{\theta}}{\theta}\right)^3\mathcal{I}\{m_n = (n+3)/2\} + O(n^{-2}) \\
&= \frac{2(n - m_n - 1)}{n}\left\{\left(\frac{\tilde{\theta}}{\theta}\right)^3 - \frac{w}{4}\left[3\left(\frac{\tilde{\theta}}{\theta}\right)^4 - 2\left(\frac{\tilde{\theta}}{\theta}\right)^2 - 1\right]\right\} \\
&\quad + \frac{1}{n}\left[1 + \left(\frac{\tilde{\theta}}{\theta}\right)^2 + 2\left(\frac{\tilde{\theta}}{\theta}\right)^3\right] \\
&\quad + \frac{1}{n}\left(\frac{\tilde{\theta}}{\theta}\right)^3\mathcal{I}\{m_n = (n+3)/2\} + O(n^{-2}),
\end{aligned}
$$

as $n \to \infty$ uniformly over $\beta_0 \le \theta, \tilde{\theta} \le \beta_1$. This proves Lemma A.5. □

LEMMA A.6. *Let $\theta > 0$ and $R_{\theta,n}$ be as in Section 2. Then as $n \to \infty$,*

$$
\frac{1}{n}\operatorname{tr}\left[\left(\frac{\partial}{\partial\theta}R_{\theta,n}^{-1}\right)R_{\theta,n}\right] = -\frac{3}{\theta} + \frac{2(\theta+2)}{\theta n} + O(n^{-2})
$$

*and*

$$
\frac{1}{n}\operatorname{tr}\left[\left(\frac{\partial^2}{\partial\theta^2}R_{\theta,n}^{-1}\right)R_{\theta,n}\right] = \frac{12}{\theta^2} - \frac{2(4\theta+9)}{\theta^2 n} + O(n^{-2}).
$$

PROOF. Since $R_{\theta,n}^{-1}R_{\theta,n}$ equals the identity matrix, we have via differentiation

$$
\begin{aligned}
\frac{1}{n}\operatorname{tr}\left[\left(\frac{\partial}{\partial\theta}R_{\theta,n}^{-1}\right)R_{\theta,n}\right] &= -\frac{1}{n}\operatorname{tr}\left[R_{\theta,n}^{-1}\left(\frac{\partial}{\partial\theta}R_{\theta,n}\right)\right] \\
&= -\frac{1}{n}\sum_{i=1}^{n}\sum_{j=1}^{n}\left(\frac{\partial}{\partial\theta}R_{\theta,n}\right)_{i,j}(R_{\theta,n}^{-1})_{j,i}
\end{aligned} \tag{39}
$$

$$= \frac{1}{n}\sum_{i=1}^{n}\sum_{j=1}^{n}\frac{(i-j)^2 we^{-|i-j|w}}{n}(R_{\theta,n}^{-1})_{j,i}$$

and

$$\begin{aligned}\frac{1}{n}\operatorname{tr}\left[\left(\frac{\partial^2}{\partial\theta^2}R_{\theta,n}^{-1}\right)R_{\theta,n}\right] &= \frac{1}{n}\sum_{i=1}^{n}\left[\left(\frac{\partial^2}{\partial\theta^2}R_{\theta,n}^{-1}\right)R_{\theta,n}\right]_{i,i}\\ &= \frac{2}{n}\sum_{i=1}^{n}\left[R_{\theta,n}^{-1}\left(\frac{\partial}{\partial\theta}R_{\theta,n}\right)R_{\theta,n}^{-1}\left(\frac{\partial}{\partial\theta}R_{\theta,n}\right)\right]_{i,i}\\ &\quad - \frac{1}{n}\sum_{i=1}^{n}\left[R_{\theta,n}^{-1}\left(\frac{\partial^2}{\partial\theta^2}R_{\theta,n}\right)\right]_{i,i},\end{aligned} \tag{40}$$

where, for all $1 \le i,j \le n$,

$$\begin{aligned}\left(\frac{\partial}{\partial\theta}R_{\theta,n}\right)_{i,j} &= -\frac{(i-j)^2 we^{-|i-j|w}}{n},\\ \left(\frac{\partial^2}{\partial\theta^2}R_{\theta,n}\right)_{i,j} &= -\frac{(i-j)^2(1-w|i-j|)e^{-|i-j|w}}{n^2}.\end{aligned} \tag{41}$$

Using (41) and Theorems 6 and 7, Lemma A.6 is proved by expanding the right-hand side of (39) and (40) as a power series in $w$ using *Mathematica*. We refer the reader to [5] for more details. □

LEMMA A.7. *Let $0 < \beta_0 \le \theta, \tilde{\theta} \le \beta_1 < \infty$. Then*

$$\frac{1}{n}\operatorname{tr}[(R_{\theta,n}^{-1}R_{\tilde{\theta},n})^2] = O(1),$$

*as $n \to \infty$ uniformly over $\beta_0 \le \theta, \tilde{\theta} \le \beta_1$.*

REMARK. The proof of Lemma A.7, though conceptually simple, is extremely tedious and the symbolic computation software *Mathematica* features significantly in the evaluation of the error terms. A detailed proof can be found in [5]. A much abbreviated proof is given below.

PROOF OF LEMMA A.7. First, it is convenient to note that

$$\begin{aligned}(R_{\theta,n}^{-1}R_{\tilde{\theta},n})_{i,j} &= (R_{\tilde{\theta},n}R_{\theta,n}^{-1})_{j,i},\\ (R_{\tilde{\theta},n}R_{\theta,n}^{-1})_{i,j} &= (R_{\tilde{\theta},n}R_{\theta,n}^{-1})_{n-i+1,n-j+1},\end{aligned}$$

whenever $1 \le i,j \le n$ and

$$[(R_{\theta,n}^{-1}R_{\tilde{\theta},n})^2]_{i,i} = [(R_{\theta,n}^{-1}R_{\tilde{\theta},n})^2]_{n-i+1,n-i+1} \qquad \forall 1 \le i \le n.$$

Also, if $n$ is an odd integer, we have

$$[(R_{\theta,n}^{-1}R_{\tilde{\theta},n})^2]_{(n+1)/2,(n+1)/2} = \left(\frac{\tilde{\theta}}{\theta}\right)^6 + O(n^{-1}),$$

as $n \to \infty$ uniformly over $\beta_0 \le \theta, \tilde{\theta} \le \beta_1$. Hence, writing $m_n = \lfloor (n+3)/2 \rfloor$, we obtain

$$\begin{aligned}
&\frac{1}{n}\operatorname{tr}[(R_{\theta,n}^{-1}R_{\tilde{\theta},n})^2] \\
&\qquad = \frac{2}{n}\sum_{i=m_n}^{n}[(R_{\theta,n}^{-1}R_{\tilde{\theta},n})^2]_{i,i} + \frac{1}{n}\left(\frac{\tilde{\theta}}{\theta}\right)^6 \mathcal{I}\{m_n = (n+3)/2\} + O(n^{-2}),
\end{aligned}$$

as $n \to \infty$ uniformly over $\beta_0 \le \theta, \tilde{\theta} \le \beta_1$. Consequently, it follows from

$$[(R_{\theta,n}^{-1}R_{\tilde{\theta},n})^2]_{j,j} = O(n) \qquad \forall\, j = n-1, n,$$

$$\frac{2}{n}\sum_{i=m_n}^{n-2}[(R_{\theta,n}^{-1}R_{\tilde{\theta},n})^2]_{i,i} = O(1),$$

that

$$\frac{1}{n}\operatorname{tr}[(R_{\theta,n}^{-1}R_{\tilde{\theta},n})^2] = O(1),$$

as $n \to \infty$ uniformly over $\beta_0 \le \theta, \tilde{\theta} \le \beta_1$. $\square$

## APPENDIX B

For $i = 1, 2$, we define

$$\begin{aligned}
C_{1,i} = {} & \frac{\hat{\tau}_{-2}\tilde{\tau}_0[w(2-6u^2+4u^4)+(1-u^2)^2]}{\tau_{-1}u^2[(b+\tau_{-1}u)\alpha_i - \tau_{-1}u]} \\
& + \frac{w\hat{\tau}_{-2}\tilde{\tau}_0(1-u^2)^2}{u[(b+\tau_{-1}u)\alpha_i - \tau_{-1}u]^2} + \frac{\hat{\tau}_0\tilde{\tau}_0 - \hat{\tau}_1\tilde{\tau}_{-1}}{\tau_{-1}u^4} \\
& - \frac{\hat{\tau}_{-2}\tilde{\tau}_{-1}[w(3-8u^2+5u^4)+(1-u^2)^2]}{\tau_{-1}u[(b+\tau_{-1}u)\alpha_i - \tau_{-1}u]} - \frac{w\hat{\tau}_{-2}\tilde{\tau}_{-1}(1-u^2)^2}{[(b+\tau_{-1}u)\alpha_i - \tau_{-1}u]^2}, \\
C_{2,i} = {} & \frac{w\hat{\tau}_{-2}(2-3u^2+u^6)[w(2-6u^2+4u^4)+(1-u^2)^2]}{2\tau_{-1}u^2[(b+\tau_{-1}u)\alpha_i - \tau_{-1}u]} \\
& + \frac{w^2\hat{\tau}_{-2}(2-3u^2+u^6)(1-u^2)^2}{2u[(b+\tau_{-1}u)\alpha_i - \tau_{-1}u]^2} + \frac{w\hat{\tau}_0(2-3u^2+u^6)}{2\tau_{-1}u^4} \\
& - \frac{w(1-u^2)^2\hat{\tau}_{-2}\tilde{\tau}_{-1}}{\tau_{-1}u[(b+\tau_{-1}u)\alpha_i - \tau_{-1}u]} - \frac{w(1-3u^4+2u^6)\tilde{\tau}_{-1}}{\tau_{-1}u^3},
\end{aligned}$$

$$C_{3,i} = \hat{\tau}_{-1} + \frac{\hat{\tau}_{-2}u[w(1-4u^2+3u^4)+(1-u^2)^2]}{(b+\tau_{-1}u)\alpha_i - \tau_{-1}u} + \frac{w\tau_{-1}u^2\hat{\tau}_{-2}(1-u^2)^2}{[(b+\tau_{-1}u)\alpha_i - \tau_{-1}u]^2},$$

$$\begin{aligned} C_{4,i} = {} & \frac{\hat{\tau}_{-1}\tilde{\tau}_0}{\tau_{-1}u^4} + \frac{\hat{\tau}_{-2}\tilde{\tau}_0[w(1-4u^2+3u^4)+(1-u^2)^2]}{\tau_{-1}u^3[(b+\tau_{-1}u)\alpha_i - \tau_{-1}u]} \\ & + \frac{w\tau_{-1}\hat{\tau}_{-2}\tilde{\tau}_0(1-u^2)^2}{\tau_{-1}u^2[(b+\tau_{-1}u)\alpha_i - \tau_{-1}u]^2} \\ & - \frac{\hat{\tau}_1\tilde{\tau}_{-2}}{\tau_{-1}u^4} - \frac{\hat{\tau}_{-2}\tilde{\tau}_{-2}[w(3-8u^2+5u^4)+(1-u^2)^2]}{\tau_{-1}u[(b+\tau_{-1}u)\alpha_i - \tau_{-1}u]} \\ & - \frac{w\hat{\tau}_{-2}\tilde{\tau}_{-2}(1-u^2)^2}{[(b+\tau_{-1}u)\alpha_i - \tau_{-1}u]^2}, \end{aligned}$$

$$\begin{aligned} C_{5,i} = {} & \frac{w\hat{\tau}_{-1}(2-3u^2+u^6)}{2\tau_{-1}u^4} \\ & + \frac{w\hat{\tau}_{-2}(2-3u^2+u^6)[w(1-4u^2+3u^4)+(1-u^2)^2]}{2\tau_{-1}u^3[(b+\tau_{-1}u)\alpha_i - \tau_{-1}u]} \\ & + \frac{w^2\hat{\tau}_{-2}(2-3u^2+u^6)(1-u^2)^2}{2u^2[(b+\tau_{-1}u)\alpha_i - \tau_{-1}u]^2} \\ & - \frac{w\hat{\tau}_{-2}\tilde{\tau}_{-2}(1-u^2)^2}{\tau_{-1}u[(b+\tau_{-1}u)\alpha_i - \tau_{-1}u]} - \frac{w(1-3u^4+2u^6)\tilde{\tau}_{-2}}{\tau_{-1}u^3}, \end{aligned}$$

$$\begin{aligned} C_{6,i} = {} & \hat{\tau}_{-1}\tilde{\tau}_{-1} - \hat{\tau}_0\tilde{\tau}_{-2} + \frac{\hat{\tau}_{-2}u\tilde{\tau}_{-1}[w(1-4u^2+3u^4)+(1-u^2)^2]}{(b+\tau_{-1}u)\alpha_i - \tau_{-1}u} \\ & - \frac{\hat{\tau}_{-2}\tilde{\tau}_{-2}u^2[w(2-6u^2+4u^4)+(1-u^2)^2]}{(b+\tau_{-1}u)\alpha_i - \tau_{-1}u} \\ & + \frac{w\tau_{-1}u^2\hat{\tau}_{-2}(\tilde{\tau}_{-1}-u\tilde{\tau}_{-2})(1-u^2)^2}{[(b+\tau_{-1}u)\alpha_i - \tau_{-1}u]^2}, \end{aligned}$$

$$C_{7,i} = \frac{\hat{\tau}_{-2}u^2[w(2-6u^2+4u^4)+(1-u^2)^2]}{[(b+\tau_{-1}u)\alpha_i - \tau_{-1}u]} + \frac{w\tau_{-1}\hat{\tau}_{-2}u^3(1-u^2)^2}{[(b+\tau_{-1}u)\alpha_i - \tau_{-1}u]^2} + \hat{\tau}_0.$$

**Acknowledgments.** This manuscript was completed when I was visiting the Department of Statistics, University of Michigan at Ann Arbor, during the fall of 2002. I would like to thank the Department, especially Professors Jeff Wu and Vijay Nair, for their very kind hospitality.

I would also like to thank Professor Michael Stein for the encouragement and the many discussions we had on fixed-domain asymptotics when I visited the University of Chicago in spring 2003.

Finally, I am very grateful to Professors Jon Wellner, Morris Eaton, an Associate Editor and two referees for their comments and suggestions that resulted in a much improved paper.

Department of Statistics
and Applied Probability
National University of Singapore
Singapore 117546
Republic of Singapore
e-mail: stalohwl@nus.edu.sg